\documentclass[11pt]{article}
\usepackage{latexsym}
\usepackage{amsfonts,amssymb,amsmath,mathtools}
\usepackage{bbold}           

\newtheorem{thm}{Theorem}[section]

\newtheorem{lem}[thm]{Lemma}
\newtheorem{pro}[thm]{Proposition}
\newtheorem{defn}[thm]{Definition}

\newcommand{\ov }{\overline }
\newcommand{\minus}{\smallsetminus}
\newcommand{\x}{\hspace{-0.025in}\times\hspace{-0.025in}}
\newcommand{\e}{\varepsilon}

\title{Evaluation problems for the Thompson group and the 
    Brin-Thompson group, and their relation to the word problem}

\author{J.C.\ Birget}

\date{\scriptsize{
16 xi 2021}}

\begin{document}
\maketitle

\begin{abstract}
The Thompson group $V$, as well as the Brin-Thompson group $2V$, is finitely
generated and can be defined as a monoid acting on bitstrings, respectively 
pairs of bitstrings. Therefore evaluation problems can be defined for $V$
and $2V$.
We show that these evaluation problems reduce to the corresponding word 
problems, and that in general, these evaluation problems are actually 
equivalent to the word problems. 
The long-input version of the evaluation problem is deterministic 
context-free and reverse deterministic context-free for $V,$ and  
{\sf P}-complete for $2V$.
\end{abstract}

\section{Introduction}

Informally, an {\em evaluation function} is of the form

\smallskip

 \ \ \   \ \ \  $E: \ (p,x) \longmapsto E_p(x)$,

\smallskip

\noindent where $p$ is a ``program'' that describes a function $E_p$, 
$x$ is a data input for $E_p$, and $E_p(x)$ is the corresponding 
data output. Here we assume that $p$, $x$, and $E_p(x)$ are strings 
over some (possibly different) alphabets.
 
The word ``input'' is ambiguous here, as $E$ has input $(p,x)$, and $E_p$ 
has input $x$. For clarity we call $x$ the {\em data input}. 

In this paper, ``function'' means {\em partial function} (unless we 
explicitly say {\em total function}, for a given source set).
So $E_p(x)$ can be undefined for some $(p,x)$.

\medskip

\noindent The {\em evaluation decision problem} of $E$ is defined as
follows.

\noindent {\sc Input:} \ $(p, x, y)$.

\noindent {\sc Question:} \ $E_p(x) = y$ ?

\smallskip

\noindent When $E_p(x)$ is undefined then $E_p(x) \ne y$; so the decision 
problem always has a {\sc yes/no} answer.
From now on we will call the evaluation decision problem simply the
{\em evaluation problem}.

\smallskip

Evaluation functions show up in many situations, e.g., in relation with 
universal Turing machines, with interpreters of programming languages, 
and more generally with exponential objects in a category. 
But here, the most relevant example here is the following.

\bigskip

\noindent {\bf Circuits:}  

\smallskip

\noindent The evaluation problem for acyclic boolean circuits, called the
{\em circuit-value problem}, is defined as follows.

\noindent {\sc Input:} \ $(C, x)$, where $C$ is (an encoding of) an 
acyclic boolean circuit with just one output wire; and $x$ is a bit-string. 

\noindent {\sc Question:} \ $C(x) = 1$ ?

\medskip

\noindent The output $C(x)$ can be undefined (if $x$ has the wrong input 
length for $C$); then the answer is {\sc no}.  
Ladner \cite{Ladner} proved that the circuit-value problem is 
{\sf P}-complete; for details about the circuit-value problem see 
\cite{Ladner, Papadim}.
We will use a more general form of the circuit-value problem, namely with 
input $(C, x, y)$, and question ``$C(x) = y$?'', where $x, y \in \{0,1\}^*$.

We can compare the circuit-value problem with the {\em circuit-equivalence 
problem}, where the input consists of two acyclic boolean circuits 
$C_1, C_2$; the question is whether $C_1$ and $C_2$ have same input-output 
function. The circuit-equivalence problem is {\sf coNP}-complete; this 
follows from the {\sf NP}-completeness of the satisfiability problem 
for boolean formulas \cite{Cook}; see \cite{Papadim} for details.

It is easy to prove that the circuit-value problem reduces to the 
circuit-equivalence problem by a many-one log-space reduction; see 
e.g.\ Section 7.1 below. But the circuit-equivalence problem does not 
reduce to the circuit-value problem, unless $\,{\sf P} = {\sf NP}$.

\bigskip

\noindent {\bf Languages and complexity:}

\smallskip

\noindent By an {\em alphabet} we mean a finite set, and we only use
alphabets that are subsets of a fixed infinite countable set.
The set of all finite strings over an alphabet $A$ is denoted by $A^*$; this
includes the empty string $\e$; we denote $A^* \minus \{\e\}$ by $A^+$.
Strings can be concatenated, which makes $A^*$ the free monoid overe $A$.  
For $m \in {\mathbb N}$, the set of strings in $A^*$ of length $m$ is denoted
by $A^m$, and the set of strings of length $\le m$ is denoted by $A^{\le m}$;
in particular, $A^0 = \{\e\}$. The length of a string $x$ is denoted by 
$|x|$. For a set $S$, the cardinal of $S$ is denoted by $|S|$. 

We will be interested in the complexity of some evaluation problems, and
we will use the well-known complexity classes and reductions below; see 
e.g.\ \cite{HU, Papadim, GHR, HemaOgi, DuKo}; for context-free languages, 
see especially \cite{Harrison, HU}. 
Since all our alphabets are subsets of a fixed infinite countable set, 
the set of Turing machines is countable, and each complexity ``class'' 
below is a countable set.

\smallskip

\noindent $\bullet$ \ \ {\sf CF} \ -- \ the {\em context-free languages}.

\smallskip

\noindent $\bullet$ \ \ {\sf coCF}  \ -- \ the {\em co-context-free 
languages}; \ {\sf coCF} $=$ 
$\,\{L \subseteq A^* : A$ is an alphabet, $A^* \minus L \in {\sf CF}\}$.

\smallskip

\noindent $\bullet$ \ \ {\sf DCF} \ -- \ the {\em deterministic context-free} 
languages.

\smallskip

\noindent $\bullet$ \ \ ${\sf DCF}^{\rm rev}$ \ -- \ the {\em reverse 
deterministic context-free} languages; \ ${\sf DCF}^{\rm rev} = \{L :$
$L^{\rm rev} \in {\sf DCF}\}$.

\smallskip

\noindent $\bullet$ \ \ {\sf log-space} computable total functions and 
languages accepted in logarithmic space.

\smallskip

\noindent $\bullet$ \ \ {\sf P} \ -- \ the set of languages accepted by 
deterministic polynomial-time Turing machines. 

\smallskip

\noindent $\bullet$ \ \  {\sf NP} \ -- \ the set of languages accepted by 
nondeterministic polynomial-time Turing machines with existential 
acceptance.

\smallskip

\noindent $\bullet$ \ \ {\sf coNP} \ -- \ the set of languages accepted by 
nondeterministic polynomial-time Turing machines with universal 
acceptance; equivalently, $\,{\sf coNP} \,=\, \{L \subseteq A^* : A$ is a 
finite alphabet, $A^* \minus L \in {\sf NP} \}$.
The set {\sf coNP} has the following useful characterization. 
For any $L \subseteq A^*$ we have: 

\smallskip

\hspace{0.2in} $L \in {\sf coNP}$ \ \ iff 

\smallskip

\hspace{0.2in} there exists a two-variable predicate $P_L(.,.)$  $\subseteq$ 
$A^* \x B^*\,$ that is decidable in deterministic 

\hspace{0.2in} polynomial-time (where $B$ is an alphabet), and 
there exists a polynomial $\pi_L(.)$, such that

\hspace{0.2in} $L \,=\, \{ x \in A^* : \ $ 
$(\forall y \in A^{\le \pi_L(|x|)}) \,P_L(x,y) \,\}$

\smallskip

\noindent For {\sf NP} a similar characterization applies, but with 
$\forall$ replaced by $\exists$.

\medskip 

\noindent To define {\em completeness} in a complexity class we use various 
reductions. Let $L_1 \subseteq A^*$ and $L_2 \subseteq B^*$ two languages. 

\smallskip

\noindent $\bullet$ \ \ A {\em many-one log-space reduction} from 
$L_1 \subseteq A^*$ to $L_2 \subseteq B^*$ is a log-space computable total 
function $f: A^* \to B^*$ such that $L_1 = f^{-1}(L_2)$; equivalently, for 
all $x \in A^*$: \ $x \in L_1 \ $ iff $ \ f(x) \in L_2$.
 
\smallskip

\noindent $\bullet$ \ \ An {\em $N$-ary conjunctive log-space reduction} 
from $L_1 \subseteq A^*$ to $L_2 \subseteq B^*$ (for some $N > 0$) is a 
log-space computable total function $\,f: x \in A^*$  $\,\longmapsto\,$ 
$f(x) = \big(f(x)_1, \,\ldots\,, f(x)_N \big)$   $\in$ 
{\Large \sf X}$_{i=1}^{^N} B^*\,$ such that for all $x \in A^*$:
 \ $x \in L_1 \ $ iff $ \ f(x)_1 \in L_2\,$ {\sc and} $ \ \ldots \ $ 
{\sc and}  $\,f(x)_N \in L_2$.

\smallskip

\noindent $\bullet$ \ \ A {\em conjunctive log-space reduction of polynomial
arity} from $L_1 \subseteq A^*$ to $L_2 \subseteq B^*$ consists of a 
polynomial $\pi(.)$ and a log-space computable total function 
$\,f: x \in A^* \,\longmapsto\,$ 
$f(x) = \big(f(x)_1, \,\ldots\,, f(x)_{\pi(|x|)} \big)$  $\in\,$ 
{\Large \sf X}$_{i=1}^{^{\pi(|x|)}} B^*\,$ such that for all $x \in A^*$:
\ $x \in L_1 \ $ iff \ (for all $i = 1, 2, \,\ldots\,, \pi(|x|))$: 
$\, f(x)_i \in L_2$.

\medskip

The complexity classes {\sf P}, {\sf NP}, and {\sf coNP} are downward closed 
under these reductions; i.e., if $L_2$ is in the class, and $L_1$ reduces to 
$L_2$, then $L_1$ is in the class.

A language $L$ is {\em complete} in a class $\cal C$ for a certain 
type of reduction iff $L \in {\cal C}$, and every language in $\cal C$ 
reduces to $L$ for this type of reduction.  

\bigskip

Since the Thompson group $V$ and the Brin-Thompson group $2V$ are finitely 
generated, and are transformation groups (acting on $\{0,1\}^{\omega}$,
respectively $2\,\{0,1\}^{\omega}$), we can also consider evaluation 
functions and evaluation problems for $V$ and $2V$. 
For the program input, a string of generators of $V$ or $2V$ is used.
For the data input, however, there is a complication: $V$ and $2V$ do not 
act (as transformation groups) on finite strings. 
Nevertheless, $V$ and $2V$ can also be defined by partial transformations 
on $\{0,1\}^*$, respectively $2\,\{0,1\}^*$, as described below.
Hence evaluation problems for $V$ and $2V$ can be defined with a string, 
respectively a pair of strings, as data input.

Regarding a generating set $\Gamma$ of the groups $V$ and $2V$ we make the 
convention that $\Gamma$ is {\em closed under inverse}; i.e., by $\Gamma^*$
we always mean $(\Gamma^{\pm 1})^*$.


\section{Evaluation problems for the Thompson group} 

\smallskip

\noindent For the definition of the Thompson group $V$ we follow 
\cite[Sect.\ 2.1]{BinG}; we will repeat some of the definitions, but not all.

For $x, p \in A^*$ we say that $p$ is a {\it prefix} of $x$ iff
$(\exists u \in A^*) \, x = pu$; this is denoted by $p \le_{\rm pref} x$.
A {\it prefix code} is any subset $P \subseteq A^*$ such that no element of 
$P$ is a prefix of another element of $P$. A {\em maximal prefix code} in 
$A^*$ is a prefix code that not a strict subset of any prefix code in $A^*$.
A {\it right ideal} is a subset $R \subseteq A^*$ such that $R A^* = R$.  
A right ideal is $R$ {\em essential} in $A^*$ iff $R$ has a non-empty 
intersection with every non-$\varnothing$ right ideal of $A^*$. For every 
right ideal $R$ there exists a unique prefix code $P$ such that $R = PA^*$; 
and $R$ is essential iff $P$ is maximal (see e.g.\ 
\cite[Lemma 8.1]{BiThomps}). 

A {\it right ideal morphism} of $A^*$ is a function $f: A^* \to A^*$, with 
domain ${\rm Dom}(f)$ and image set ${\rm Im}(f)$, such that for all 
$x \in {\rm Dom}(f)$ and all $w \in A^*$: \ $f(xw) = f(x) \ w$. 
The unique prefix code that generates the right ideal ${\rm Dom}(f)$ is 
denoted by ${\rm domC}(f)$, and is called the {\it domain code} of $f$; the 
unique prefix code that generates the right ideal ${\rm Im}(f)$ is denoted by 
${\rm imC}(f)$, and is called the {\it image code}. In order to define $V$ 
we first introduce the inverse monoid

\medskip

${\cal RI}_A^{\sf fin}$ $\, = \,$  $\{ f : f$ is a right ideal morphism of
    $A^*$ such that $f$ is {\em injective}, and

\hspace{0.9in}
    ${\rm domC}(f)$ and ${\rm imC}(f)$ are {\em finite maximal}
    prefix codes\}.

\medskip

\noindent Every $f \in {\cal RI}_A^{\sf fin}$ has a unique maximum extension
to an element of ${\cal RI}_A^{\sf fin}$ (by \cite[Prop.\ 2.1]{BiThomps}).

We define the Higman-Thompson group $G_{k,1}$ (where $k = |A|$) as follows:
As a set, $G_{k,1}$ consists of the right ideal morphisms $f$  $\in$ 
${\cal RI}_A^{\sf fin}$ that are maximum extensions in 
${\cal RI}_A^{\sf fin}$; so $G_{k,1}$  $\subseteq$ ${\cal RI}_A^{\sf fin}\,$ 
(as sets).
The multiplication in $G_{k,1}$ consists of composition, followed by maximum
extension. The Thompson group $V$ is $G_{2,1}$.

\smallskip

There are other characterizations of $G_{k,1}$; we give two more, one based 
on $\,\equiv_{\rm end}$, and one based on a faithful action on $A^{\omega}$.
For $f \in {\cal RI}_A^{\sf fin}$, $p \in {\rm domC}(f)$, and 
$u \in A^{\omega}$, we define $f(pu) = f(p) \, u$.

\smallskip

\noindent {\small (1)} The group $G_{k,1}$ is also a homomorphic image of 
${\cal RI}_A^{\sf fin}$. For $f_1, f_2 \in$  ${\cal RI}_A^{\sf fin}$ we 
define the congruence $\,\equiv_{\rm end}$ as follows: 

\smallskip

 \ \ \  \ \ \  $f_1 \equiv_{\rm end} f_2$ \ \ iff \ \ $f_1$ and $f_2$ agree 
on $\,{\rm Dom}(f_1) \cap {\rm Dom}(f_2)$. 

\medskip

\noindent Then $G_{k,1}$ is isomorphic to
$ \ {\cal RI}_A^{\sf fin}\!/\!\!\equiv_{\rm end}$.

\smallskip

\noindent {\small (2)} The group $G_{k,1}$ is isomorphic to the action monoid 
of the action of ${\cal RI}_A^{\sf fin}$ on $A^{\omega}$. Indeed,  
$\,f_1 \equiv_{\rm end} f_2$ \ iff \ the actions of $f_1$ and $f_2$ on 
$A^{\omega}$ are the same. 

\smallskip

\noindent In summary, we have three equivalent definitions of $G_{k,1}$: 

\noindent {\small (0)} As the subset (not subgroup) of 
${\cal RI}_A^{\sf fin}$ consisting of maximally extended right-ideal 
morphisms, with multiplication consisting of composition, followed by maximum 
extension. The definition in \cite{CFP} is a numerical coding of this 
definition; see \cite[Sect.\ 2.1]{BinG}.

\noindent {\small (1)} As $\,{\cal RI}_A^{\sf fin}\!/\!\!\equiv_{\rm end}$.

\noindent {\small (2)} As the action monoid of ${\cal RI}_A^{\sf fin}$ on
$A^{\omega}$.

\medskip

Every element $f \in {\cal RI}_A^{\sf fin}$ (and in particular, every
$f \in G_{k,1}$) is determined by the restriction of $f$ to ${\rm domC}(f)$.
This restriction $\,f|_{{\rm domC}(f)}: {\rm domC}(f) \to {\rm imC}(f)\,$ is
a finite bijection, called the {\it table} of $f$ (see \cite{Hig74}).
When we use tables we do not always assume that $f$ is a maximum extension;
but the maximum extension can easily be found from the table.

\bigskip

\noindent To define an {\em evaluation function} for $V$ we first choose a 
finite generating set $\Gamma_{\!1}$. The evaluation function for $V$ over 
$\Gamma_{\!1}$ is 

\medskip

 \ \ \   \ \ \   $E: \ (w,x) \in \Gamma_{\!\!1}^* \x \{0,1\}^*$ 
$ \ \longmapsto \ $  $E_w(x) \in \{0,1\}^*$,

\medskip

\noindent where, if $w = w_n \,\ldots\, w_1$ with $w_n, \ \ldots \ , w_1$ 
$\in$ $\Gamma_{\!1}$, and $E_w$ is the element of $V$ generated by the string
$w$.  By \cite[Prop.\ 2.1]{BiThomps}, $E_w$ is the maximum extension of 
$\,w_n \circ \,\ldots\, \circ w_1(.) \in {\cal RI}_A^{\sf fin}$ to a 
right-ideal morphism in ${\cal RI}_A^{\sf fin}$, and this maximal extension 
is unique.

For the Thompson group $V$ it is useful to view the data input as a {\em
stack} (push-down store); in the Brin-Thompson group, the data input is a
pair of stacks. An element of $V$ and $2V$ changes the top bits of the
stack(s).

\begin{defn} \label{DEFevVGam} {\bf (evaluation problem of $V$ over
$\Gamma_{\!1}$).}

\noindent The {\em evaluation problem} of $V$ over a finite generating set
$\Gamma_{\!1}$ is specified as follows.

\noindent {\sc Input:} \ $(w, x, y)$  $\in$ 
$\Gamma_{\!\!1}^* \x \{0,1\}^* \x \{0,1\}^*$.

\noindent {\sc Question:} \ $E_w(x) = y$ ? 

\smallskip

\noindent In other words, this problem is the set
$\,\{(w, x, y) \in \Gamma_{\!\!1}^* \x \{0,1\}^* \x \{0,1\}^*:\,$ 
$E_w(x) = y\}$.
\end{defn}

\noindent The following gives a connection between the action of $V$ on 
finite strings and on infinite one.
First a general Lemma about the relation between $\{0,1\}^*$ and
$\{0,1\}^{\omega}$.

\begin{lem} \label{LEMequalityEquivalents}
 \ The following equivalences hold for every $x, y \in \{0,1\}^*$:

\medskip

 \ \ \   \ \ \ $x = y\,$
 \ \ \ $\Leftrightarrow$
 \ \ \ $x\,\{0,1\}^* \,=\, y \,\{0,1\}^*$
 \ \ \ $\Leftrightarrow$
 \ \ \ $x \,\{0,1\}^{\omega} \,=\, y \,\{0,1\}^{\omega}$.

\medskip

\noindent Similarly,

\smallskip

 \ \ \  \ \ \ $y \le_{\rm pref} x$
 \ \ \ $\Leftrightarrow$
 \ \ \ $x \,\{0,1\}^* \,\subseteq\, y \,\{0,1\}^*$
 \ \ \ $\Leftrightarrow$
 \ \ \ $x \,\{0,1\}^{\omega} \,\subseteq\, y \,\{0,1\}^{\omega}$.
\end{lem}
{\sc Proof.} Obviously, $x = y$ implies the other equalities, and 
$y \le_{\rm pref} x$ implies the inclusions.
Conversely, if $x \ \{0,1\}^* \subseteq y \ \{0,1\}^*$, or
$x \ \{0,1\}^{\omega}$ $\subseteq$  $y \ \{0,1\}^{\omega}$, then
$x \in y \ \{0,1\}^*$, hence $y \le_{\rm pref} x$. 

Symmetrically, if we have the other inclusions then we also have 
$x \le_{\rm pref} y$, hence, $x = y$.
 \ \ \ $\Box$

\begin{lem} \label{LEMfinInfin}
 \ For every element $f \in V\,$ ($\,\subseteq$ 
${\cal RI}_{\{0,1\}}^{\sf fin}$) and every $x, y \in \{0,1\}^*$, the 
following are equivalent:

\smallskip

\noindent {\small \rm (1)} \ \ $f(x) = y$.

\smallskip

\noindent {\small \rm (2)} \ \ For all $z \in \{0,1\}^*$: \ $f(xz) = yz$.

\smallskip

\noindent {\small \rm (3)} \ \ There exists a maximal prefix code 
$P \subseteq \{0,1\}^*$ such that for all $z \in P$: \ $f(xz) = yz$. 

\smallskip

\noindent {\small \rm (4)} \ \ For all $u \in \{0,1\}^{\omega}$: 
 \ $f(x u) = y u$.
\end{lem}
{\sc Proof.} The implications (1) $\Rightarrow$ (2), (1) $\Rightarrow$ (3),
and (1) $\Rightarrow$ (4) are obvious. And (1) $\Leftarrow$ (2) is obtained 
by taking $z = \e$.

(1) $\Leftarrow$ (3): Suppose $f(xz) = yz\,$ for all $z \in P$, where $P$ is 
a maximal prefix code.  Then $f(x) = y$ holds in the maximum extension of 
$f$ to a right ideal morphism. Since $f$ is already maximally extended (as 
$f \in V$), we have $f(x) = y$.

(1) $\Leftarrow$ (4):  Assume $f(xu) = y u$.  Since ${\rm domC}(f)$ is a 
maximal finite prefix code, every $u \in \{0,1\}^{\omega}$ has a prefix 
$z_{x,u}$ such that $f(x z_{x,u})$ is defined. Let 
$v_{x,u} \in \{0,1\}^{\omega}$ be such that $u = z_{x,u} v_{x,u}$. So, 
$yu$  $=$ $f(xu)$ $=$ $f(x z_{x,u} v_{x,u})$  $=$  $f(x z_{x,u})\,v_{x,u}$ 
$=$  $y z_{x,u} v_{x,u}$. Hence, since $v_{x,u}$ can be any element of 
$\{0,1\}^{\omega}$, Lemma \ref{LEMequalityEquivalents} implies that for every
$u \in \{0,1\}^{\omega}$: \ $f(x z_{x,u}) = y z_{x,u}$.  
Since $z_{x,u}$ exists for every $x \in {\rm domC}(f)$ and every $u$, it 
follows that $\,\{z_{x,u} : x \in {\rm domC}(f), \ u \in \{0,1\}^{\omega}\}$ 
$\cdot$  $\{0,1\}^*\,$ is an essential right ideal. Let $P$ be the maximal 
prefix code that generates this essential right ideal. Then for all 
$p \in P$: $\,f(xp) = yp$.  Now (3) holds, which implies (1).
 \ \ \ $\Box$ 

\bigskip

\noindent {\bf Notation:} \\    
For $x = (x_1, \,\ldots\,, x_n) \in n\,A^*$:
 \ \ ${\rm maxlen}(x) = \max\{|x_i| : i = 1, \,\ldots\,, n\}$.

\noindent For a finite set $S \subseteq n\,A^*$:
 \ \ ${\rm maxlen}(S) = \max\{ {\rm maxlen}(x) : x \in S\}$.

\noindent For a table $F: P \to Q$:
 \ \ ${\rm maxlen}(F) = \max\{ {\rm maxlen}(x) : x \in P \cup Q\}$.

\noindent For a set of tables $S$, 
 \ \ ${\rm maxlen}(S) = \max\{ {\rm maxlen}(F) : F \in S\}$.

\bigskip 

\noindent {\bf Long versus short data inputs:} When $V$ is defined by
functions on $\{0,1\}^*$, $E_w$ is an injective function; but $E_w(x)$ is 
not defined for all $x \in \{0,1\}^*$, except when $E_w$ is the identity 
function. Moreover, $E_w$ can be a strict extension of 
$w_n \circ \,\ldots\, \circ w_1(x)$.

If $x$ is long enough then $\,E_w(x) = w_n \circ \,\ldots\, \circ w_1(x)$;
i.e., $E_w(x)$ is obtained by simply applying the generators 
$w_i \in \Gamma_{\!1}$ to the data input, one generator after another. 
A sufficient (but not necessary) condition for this is given by the following
\cite[Cor.\ 3.7]{BiThomps}: Let 

\smallskip

 \ \ \  \ \ \ $c_{_{\Gamma_1}} = {\rm maxlen}(\Gamma_{\!1})$, 

\smallskip

\noindent i.e., the length of the longest bit-string in the tables of the 
generators in $\Gamma_{\!1}$. Then for every $w \in \Gamma_{\!\!1}^*$ and
$x \in \{0,1\}^*$:

\medskip

 \ \ \  \ \ \ if $ \ |x| \,\ge\, c_{_{\Gamma_1}} \,|w|$ \ then 
 \ $w_n \circ \,\ldots\, \circ w_1(x) \ $ is defined;

\medskip

\noindent and in that case, $\,E_w(x) = w_n \circ \,\ldots\, \circ w_1(x)$.

\begin{defn} \label{DEFlongshortdata} {\bf (long versus short data inputs).}

\smallskip

\noindent Let $w \in \Gamma_{\!\!1}^*$ with $w = w_n \ldots w_1\,$ and 
$\,w_n, \,\ldots\,, w_1$ $\in$ $\Gamma_{\!1}$.  

We call a data input $x \in \{0,1\}^*$ a {\em long data input for} the 
word $w$
 \ iff \ $w_n \circ \,\ldots\, \circ w_1(x)\,$ is defined (and then
$\,E_w(x) = w_n \circ \,\ldots\, \circ w_1(x)$).

A data input $x \in {\rm Dom}(E_w)$ that is not long is called
a {\em short data input for} the word $w$.
\end{defn}
In summary, for $w \in \Gamma_{\!\!1}^*$ there are three kinds of data 
inputs $x$: 

\noindent (1) data inputs that are {\em too short}, i.e., 
$x \not\in {\rm Dom}(E_w)$;

\noindent (2) data inputs that are {\em short}, i.e., 
$x \in {\rm Dom}(E_w)$, but 
$w_n \circ \,\ldots\, \circ w_1(x)$ is undefined.

\noindent (3) data inputs that are {\em long}, i.e., 
$x \in {\rm Dom}(w_n \circ \,\ldots\, \circ w_1(.))$.

\bigskip

For the Thompson group $V$, one can also consider a {\em circuit-like
generating} set $\Gamma_{\!1} \cup \tau$, and define the evaluation function,
evaluation problem, and domain membership problem of $V$ over
$\Gamma_{\!1} \cup \tau$. 
Here, $\Gamma_{\!1}$ is any finite generating set of $V$, and $\tau$ is the 
set of {\em bit-transpositions}, i.e., $\tau = \{\tau_{i,i+1} : i \ge 1\}$, 
where $\tau_{i,i+1}$ is the right ideal morphism of $\{0,1\}^*$ defined by  
$\,\tau_{i,i+1}(u x_i x_{i+1} v)$ $=$ $u x_{i+1} x_i v\,$ for all 
$u \in \{0,1\}^{i-1},$ $v \in \{0,1\}^*,$ and $x_i, x_{i+1} \in \{0,1\}$; 
$\tau(z)$ is undefined if $|z| \le i$. See \cite{BiCoNP,BinG} for the 
connection between $\tau$ and acyclic boolean circuits.

\begin{defn} \label{DEFevVGamTau} {\bf (evaluation problem of $V$ over
$\Gamma_{\!1} \cup \tau$).}

\noindent The {\em evaluation problem} of $V$ over $\Gamma_{\!1} \cup \tau$ 
is specified as follows.

\noindent {\sc Input:} \ $(w, x, y)$  $\in$
$(\Gamma_{\!\!1} \cup \tau)^* \x \{0,1\}^* \x \{0,1\}^*$.

\noindent {\sc Question:} \ $E_w(x) = y$ ?
\end{defn}
When we consider the evaluation problem and the word problem of $V$ over
$\Gamma_{\!1} \cup \tau$, the elements of $\tau$ are encoded over a finite
alphabet. We will use the following: $\,\tau_{j,j+1}$ is encoded by
$a b^{j+1} a$. We assume that $\{a,b\}$, $\{0,1\}$ and
$\Gamma_{\!1}$ have empty intersection two by two.

\bigskip

{\em Long} and {\em short} data inputs over $\Gamma_{\!\!1} \cup \tau$ are 
defined in the same way as over $\Gamma_{\!1}$.
For $w \in (\Gamma_{\!1} \cup \tau)^*$, a sufficient condition for $x \in$ 
$\{0,1\}^*$ to be a long data input for $w$ is 

\medskip

 \ \ \  \ \ \  $|x| \,\ge\, c_{_{\Gamma_1, w}} \, |w|$;  

\smallskip

\noindent here 

 \ \ \  \ \ \  $c_{_{\Gamma_1, w}}$  $=$ 
$\max \{ c_{_{\Gamma_1}},\, {\rm maxindex}_{\tau}(w) \}$, 

\smallskip

\noindent where $c_{_{\Gamma_1}}$ is as above, and 

\smallskip

 \ \ \  \ \ \  ${\rm maxindex}_{\tau}(w) = $ 
$\max \{i \in {\mathbb N}_{>0} : \tau_{i-1, i}$ occurs in $w \}\,$ 

\smallskip

\noindent (i.e., the largest subscript of any element of $\tau$ that 
occurs in $w$).

\begin{defn} \label{DEFevVGamLong}
 \ The evaluation problem for $V$ over $\Gamma_{\!1}$ (or over
$\Gamma_{\!1} \cup \tau$) for {\em long data inputs} is defined as follows.

\noindent {\sc Input:} \ $(w,x,y)$  $\in$
$\Gamma_{\!\!1}^* \x \{0,1\}^* \x \{0,1\}^*$ 
 \ \big(or $(\Gamma_{\!\!1} \cup \tau)^* \x \{0,1\}^* \x \{0,1\}^*$\big), 
where $w = w_n\,\ldots\,w_1$, with $w_n, \,\ldots\,, w_1 \in \Gamma_{\!1}$ 
\big(or $\Gamma_{\!\!1} \cup \tau$\big).

\noindent {\sc Question:} \ $w_n \circ \,\ldots\, \circ w_1(x) = y$ ?
\end{defn}
Remark: The question ``$E_w(x) = y$?'' is equivalent to the question
``$w_n \circ \,\ldots\, \circ w_1(x) = y$?'' iff $x$ is long for $w$. 
An answer {\sc yes} to the question 
``$w_n \circ \,\ldots\, \circ w_1(x) = y$?''
implies that $x$ is long for $w$. 
The evaluation problem is not equivalent to the evaluation problem for long 
data inputs in general.

\section{Evaluation problems for the Brin-Thompson group} 

For the definition of the Brin-Thompson group $2V$ we follow 
Sections 2.2 and 2.3 and Def.\ 2.28 in \cite{BinG}, but we will not repeat 
everything.

The $n$-fold cartesian products \ {\large \sf X}$_{_{i=1}}^{^n} A^*\,$ 
and \ {\large \sf X}$_{_{i=1}}^{^n} A^{\omega}\,$ are denoted by $nA^*$,
respectively $nA^{\omega}$.
Multiplication in $\, nA^*$ is done coordinatewise, i.e., $nA^*$ is the
direct product of $n$ copies of the free monoid $A^*$. For $u \in nA^*$ we
denote the coordinates of $u$ by $u_i \in A^*$ for $1 \le i \le n$; i.e.,
$u = (u_1, \, \ldots, u_n)$.
The {\em initial factor order} on $nA^*$ is defined as follows for 
$u, v \in nA^*$: \ $u \le_{\rm init} v$ \ iff \ there exists $x \in nA^*$ 
such that $u x = v$.  Clearly, $u \le_{\rm init} v$ iff
$u_i \le_{\rm pref} v_i$ for all $i= 1, \, \ldots, n$.

In $n A^*$, similarly to $A^*$, we have the concepts of right ideal,
essential right ideal, and generating set of a right ideal.
An {\em initial factor code} is a set $S \subseteq nA^*$ such that no element
of $S$ is an initial factor of another element of $S$.  Every right ideal is 
generated by a unique initial factor code \cite[Lemma 2.7(1)]{BinG}.

An {\em essential initial factor code} is, by definition, an initial factor 
code $S$ such that $S\, nA^*$ is an essential right ideal. It is easy to see
that every maximal initial factor code is essential. The converse does not 
hold. E.g., for $A = \{0,1\}$, $\,S = \{(\e,0),\, (0,\e),\, (1,1)\}\,$ 
is essential, but not maximal. Indeed, let $P \subseteq \{0,1\}^*$ be any 
every finite maximal prefix code with $P \ne \{\e\}$; let 
$1P = \{1p : p \in P\}$.  Then $ \ S \,\cup\, (1P \x \{\e\})\,$ is an 
essential initial factor code that has $S$ as a strict subset.

\smallskip

The {\em join} $\,u \vee v\,$ of $u,v \in nA^*$ is defined to be the unique
$\,\le_{\rm init}$-minimum common upper bound of $u$ and $v$.
Of course, $u \vee v$ does not always exist.
The join $\,u \vee v\,$ of $u = (u_1, \, \ldots, u_n)$ and 
$v = (v_1, \, \ldots, v_n)$  $\in nA^*$ is characterized by the following 
equivalent statements (where $\|_{\rm pref}$ denotes prefix comparability):

\noindent
{\rm (1)} \ \ $u \vee v \,$  exists; \\
{\rm (2)} \ \ $u$ and $v$ have a common upper bound for $\, \le_{\rm init}$,
 \ i.e., $ \, (\exists z) \, [\, u \le_{\rm init} z$  $ \,{\rm and}\, $
$v \le_{\rm init} z \,]$; \\
{\rm (3)} \ \ for all $i = 1, \, \ldots, n$:
 \ $u_i \,\|_{\rm pref}\, v_i \ $ in $A^*$.

\noindent Moreover, if
$\, u \vee v = ((u \vee v)_i : i =   1, \, \ldots, n) \,$ exists,
then $\,(u \vee v)_i \,=\, u_i \ $ if $\,v_i \le_{\rm pref} u_i$,
and $\,(u \vee v)_i \,=\,v_i \ $ if $\,u_i \le_{\rm pref} v_i$; \ see
\cite[Lemma 2.5]{BinG}.

\smallskip

A set $S \subseteq nA^*$ is {\em joinless} iff no two elements of $S$
have a join with respect to $\, \le_{\rm init}$.
Joinless sets will be called {\em joinless codes}, since they are a special 
case of initial factor codes.
For $n \ge 2$, not every initial factor code is joinless.

A {\em maximal joinless code} is a joinless code $S$ such that for all 
$x \in nA^*$, $\,S \cup \{x\}\,$ is not joinless. 
A right ideal $R$ is called {\em joinless generated} iff the unique initial
factor code that generates $R$ is joinless \cite[Def.\ 2.4 and Lemma 
2.7(2)]{BinG}.
Every maximal joinless code is an essential initial factor code 
\cite[Lemma 2.9]{BinG}; but not every an essential initial factor code is
joinless; see the example above, and \cite[Remark after Lemma 2.9]{BinG}. 
Hence for joinless codes (and in particular for prefix codes of $A^*$),
essential is equivalent to maximal; but for initial factor codes in general,
we saw above that essential is not equivalent to maximal. 

We have the following fact: If $C_1$ and $C_2$ are joinless codes then the      
right ideal $\,C_1\,nA^* \,\cap\, C_2\,nA^*\,$ is generated by the joinless 
code $\,C_1 \vee C_2$  $=$ 
$\{c_1 \vee c_2 :\,$  $c_1 \in C_1,\, c_2 \in C_2\}$; moreover, 
$C_1 \vee C_2$ is a maximal joinless code iff $C_1$ and $C_2$ are both 
maximal joinless codes \ \cite[Prop.\ 2.18]{BinG}.

In summary, in $nA^*$ with $n \ge 2$ there are two different generalizations
of the concept of prefix code, namely the initial factor codes and the 
joinless codes. And there are two different generalizations of maximal prefix 
codes, namely the essential initial factor codes and the maximal joinless 
codes.  
Initial factor codes are closely related to right ideals, whereas joinless 
codes are crucial for defining right ideal morphisms.

\smallskip

Just as for $A^*,$ one defines the concepts of {\it right ideal morphism},
domain code, and image code in $nA^*$.  At first
we only consider domain and image codes that are {\it joinless}.
Indeed, if $P \subseteq nA^*$ is not joinless, the definition of right ideal
morphisms on $P$ can be inconsistent, i.e., the morphisms might not be 
functions; see Prop.\ \ref{LEMuniqext} below.

\smallskip

\noindent Before we get to $n G_{k,1}$ we define the following monoid:

\medskip

\hspace{0.2in}
$n {\cal RI}_A^{\sf fin}$  $ \ = \ $
$\{f : \, f$ is a right ideal morphism of $nA^*$ such that $f$ is injective,

\hspace{1.28in} and ${\rm domC}(f)$ and ${\rm imC}(f)$
are {\em finite, maximal, joinless} codes\} .

\medskip

\noindent Any element of $n {\cal RI}_A^{\sf fin}$ is determined by a 
bijection $F$: $P \to Q$ between finite maximal joinless codes
$P, Q \subseteq nA^*$; such a bijection is called a {\em table}.
Conversely, every table determines an element of $n {\cal RI}_A^{\sf fin}$.

Two right ideal morphisms $f, g \in n {\cal RI}_A^{\sf fin}$ are called
{\em end-equivalent} \ iff \ $f$ and $g$ agree on
$\,{\rm Dom}(f) \,\cap\, {\rm Dom}(g)$.
This will be denoted by $f \equiv_{\rm end} g$. 
By \cite[Lemma 2.24]{BinG}: $\,f \equiv_{\rm end} g\,$ iff $\,f$ and 
$g$ have the same action on $\{0,1\}^{\omega}$.

\bigskip 

\noindent {\bf Remarks on maximal extensions of right ideal morphisms 
in $n A^*$ when $n \ge 2$:}

\smallskip

\noindent (1) If one starts out with a table $F: P \to Q$, where $F$ is a 
bijection between finite initial factor codes $P$ and $Q$, then the 
extension of $F$ to a ``right ideal morphism'' will not always be a function.
E.g., for $A = \{0,1\}$, the ``right ideal morphism'' given by the table 
$\,F =$  $\{\big((0,\e), (00,\e)\big),$  $\big((\e,0),(01,\e)\big),$ 
$\big((1,1), (1,\e)\big)\}\,$ is not a function, since 
$\,F((0,0))$ $=$  $F((0,\e))\cdot(\e,0)$  $=$  $(00,\e)\cdot(\e,0) = (00,0)$,
and also $F((0,0))$  $=$  $F((\e,0))\cdot(0,\e)$  $=$ $(01,\e)\cdot((0,\e)$ 
$=$  $(010, \e)$.

\smallskip

\noindent (2) If $P$ and $Q$ are joinless codes of equal cardinality, then 
any bijection $F: P \to Q$ determines a right ideal morphism 
$f: P\,nA^* \to Q\,nA^*\,$ (which is a function); and if $P$ and $Q$ are 
maximal joinless codes then $f$ belongs to $n \, {\cal RI}_A^{\sf fin}$. 

\smallskip

\noindent (3) If $f: P\,nA^* \to Q\,nA^*\,$ belongs to 
$n\,{\cal RI}_A^{\sf fin}$, and if $P$ and $Q$ are maximal joinless codes, 
then $f$ might be extendable to a right ideal morphism whose domain and 
image codes are essential finite initial factor codes that are not joinless.

\medskip

\noindent {\sf Example} \cite[Lemma 2.27 and Fig.\ 2]{BinG}: 

\smallskip

\noindent For $2\,\{0,1\}^*$, let 

\smallskip

$F \,=\, \{\big((0,0),(0,0)\big)$, $\,\big((1,0),(1,0)\big),$
$\,\big((0,1),(0,1)\big),$  $\,\big((1,10),(1,11)\big),$ 
$\,\big((1,11),(1,10)\big)\}$.

\smallskip

\noindent Then $F$ determines a right ideal morphism in 
$2\,{\cal RI}_2^{\sf fin}$ that can be extended to 

\smallskip

$F_1 \,=\, \{\big((\e,0),(\e,0)\big),$
$\big((0,1),(0,1)\big),$  $\big((1,10),(1,11)\big),$  
$\big( (1,11),(1,10)\big)\}\,$ $\in 2\,{\cal RI}_2^{\sf fin}$, \ \ or to 

\smallskip

$F_2 \,=\, \{\big((0,\e),(0,\e)\big),$
$\big((1,0),(1,0)\big),$  $\big((1,10),(1,11)\big),$  
$\big( (1,11),(1,10)\big)\}\,$ $\in$  $\,2\,{\cal RI}_2^{\sf fin}$.

\smallskip

\noindent $F_1$ and $F_2$ have no common extension to an element of 
$\,2\,{\cal RI}_2^{\sf fin}$. But both $F_1$ and $F_2$ can be further 
extended to

\smallskip

$F_{12} \,=\, \{\big((\e,0),(\e,0)\big),$  $\big((0,\e),(0,\e)\big),$
$\big((1,10),(1,11)\big),$
$\big( (1,11),(1,10)\big)\}$. 

\smallskip

\noindent Here, $F_{12} \not\in 2\,{\cal RI}_2^{\sf fin}\,$ (its domain
and image codes are not joinless); but $F_{12}$ is nevertheless a welldefined
right ideal morphism of $2\,\{0,1\}^*$. And $F_{12}$ is the unique maximum
extension of $F$ to a right ideal morphism of $2\,\{0,1\}^*$. 
 \ \ \  \ \ \  [End, Remarks.] 

\medskip

In the above example we see that there are two different kinds of maximal 
extensions of a right ideal morphism $f$ in $n\,{\cal RI}_A^{\sf fin}$: 

(1) non-unique maximal extensions of $f$ to elements of 
$n\,{\cal RI}_A^{\sf fin}$;

(2) a unique maximal extension of $f$ to a right ideal of $nA^*\,$ (beyond
$n\,{\cal RI}_A^{\sf fin}$).  

\noindent In general, Prop.\ \ref{LEMuniqext} will give the connection 
between $n{\cal RI}_A^{\sf fin}$ and general right ideal morphisms. 

\medskip

In \cite{BinG} and \cite{BinGk1} the definition of $n{\cal RI}_A^{\sf fin}$ 
is based on maximal joinless codes, and this is appropriate when $nV$ and 
$n{\cal RI}_A^{\sf fin}$ are considered by themselves. But for the evaluation 
problem, the action of $n{\cal RI}_A^{\sf fin}$ and $nV$ on strings matters, 
and the fact that right ideal morphisms do not have unique extensions in 
$n{\cal RI}_A^{\sf fin}$ is a problem. On the other hand, every essential 
right ideal morphism has a unique extension to an essential right ideal 
morphism, generated by an essential initial factor code; so it helps now
to consider right ideal morphisms in that setting.

By Prop.\ \ref{PROfinitcodedec} we can efficiently decide whether a table
$F: P \to Q$, where $P, Q \subseteq nA^*$ are finite initial factor codes,
describes a right ideal morphism that is welldefined (i.e., a function),
or that is total, or injective, or surjective.

\begin{lem} \label{LEMpropRI} 
 \ Let $\,f: P \,nA^* \to Q\, nA^*$ be a welldefined bijective right ideal 
morphism, described by a bijective table $\,F: P \to Q$, where 
$P, Q \subseteq nA^*$ are finite {\em essential initial factor codes}.

\smallskip

\noindent {\small \rm (1)} \ Then $\,f^{-1}: Q\, nA^* \to P \,nA^*$ is also 
a bijective right ideal morphism, with table $\,F^{-1}: Q \to P$.

\smallskip

\noindent {\small \rm (2)} \ If $S \subseteq P \,nA^*$ is an (essential) 
initial factor code then $f(S)$ is also an (essential) initial factor code.

\smallskip

\noindent {\small \rm (3)} \ If $S \subseteq P \,nA^*$ is a (maximal) 
joinless code, then $f(S)$ is also a (maximal) joinless code.
\end{lem}
{\sc Proof.} (1) Since $f$ is bijective, $f^{-1}: Q\, nA^* \to P \,nA^*$
exists as a (bijective) function. Let us prove that it is a right ideal 
morphism. For any $qx \in Q\, nA^*$, with $q \in Q$ and $x \in nA^*$,
we have $f^{-1}(qx) = pz \in P \,nA^*$, for some $p \in P$ and 
$z \in nA^*$. Applying $f$ yields $qx = f(f^{-1}(qx)) = f(pz) = f(p) \,z$, 
where $f(p) \in Q$.  

Since this holds for all $x$, we have (when $x = (\e)^n$): $q = f(p) \,z'$ 
for some $z'$. Since $Q$ is an initial factor code and both $q$ and $f(p)$ 
are in $Q$, it follows that $q = f(p)$; hence $f^{-1}(q) = p$. 

Now, going back to $qx = f(p) \,z$, we use $q = f(p)$ and cancelativity of
$nA^*$ to obtain $x = z$.
Thus, $f^{-1}(qx) = pz = f^{-1}(q) \, z = f^{-1}(q) \,x$; so $f^{-1}$ is a
right ideal morphism (given by the table $F^{-1}$). 
 
\smallskip

\noindent (2) By contraposition, if $f(S)$ is not an initial factor code 
then for some $s_1, s_2 \in S$ we have
$f(s_1) <_{\rm init} f(s_2) = f(s_1) \, u\,$ (for some $u \in nA^*$
$\minus$ $\{\e\}^n$).
Then $s_2 = f^{-1}f(s_2) = f^{-1}(f(s_1) \, u) = f^{-1}(f(s_1)) \ u\,$
(the latter by (1) above). Hence $s_1 <_{\rm init} s_2 = s_1 u$, since 
$u \ne (\e)^n$. This implies that $S$ is not an initial factor code.

Since $P$ is essential we have: $S$ is essential iff for all 
$p \in P$, $w \in nA^*$: $ \ S\,nA^* \,\cap\, pw\,nA^* \ne \varnothing$;
i.e. for all all $p \in P$, $w \in nA^*$ there exist $s \in S$,
$x, y \in nA^*$ such that $s x = pw y$. Hence for all $p \in P$, 
$w \in nA^*$ there exist $f(s) \in f(S)$ and $x, y \in nA^*$ such that 
$\,f(s)\,x = f(p)\,wy$.
And $f(p)$ ranges over all of $Q$.  So for all $q \in Q$, $w \in nA^*$: 
$ \ f(S)\,nA^* \,\cap\, qw\,nA^*$  $\ne$  $\varnothing$. \ Hence (since
$Q$ is essential),  $f(S)$ is essential. 

\smallskip

\noindent (3) By contraposition, if $f(S)$ is not joinless then 
$\,f(s_1) \vee f(s_2)\,$ exists in $f(P \, nA^*)$, for some $s_1, s_2 \in S$ 
with $s_1 \ne s_2$.  Then $f(s_1) \vee f(s_2) = f(j)\,$ for some 
$j \in P \,nA^*$, so $f(j) = f(s_1) \, u_1 = f(s_2) \, u_2$ for some 
$u_1, u_2 \in nA^*$. Hence, applying $f^{-1}$ yields $j = s_1 u_1 = s_2 u_2$,
so $\,s_1 \vee s_2\,$ exists, So $S$ is not joinless.

By contraposition, if $f(S)$ is joinless but not maximal joinless then there 
exists $f(z) \in Q\,nA^*$ such that $f(S) \cup \{f(z)\}$  $=$  
$f(S \cup \{z\})$ is joinless. Then, applying $f^{-1}$ and the result from 
the paragraph above) shows that $S \cup \{z\}$ joinless. Hence, $S$ is not 
a maximal joinless code.  
 \ \ \ $\Box$

\begin{pro} \label{LEMuniqext} {\bf (maximal extensions of a right ideal
morphism).}

\smallskip

\noindent {\small \rm (1)} For $n \ge 2$ there exist right ideal morphisms 
in $n{\cal RI}_A^{\sf fin}$ that can be extended in more than one way to 
a maximal right ideal morphisms in $n{\cal RI}_A^{\sf fin}$.

\smallskip

\noindent {\small \rm (2)} Every right ideal morphism in 
$n{\cal RI}_A^{\sf fin}$ has a {\em unique} maximum extension to a right 
ideal morphism of $n\,A^*$, with domain code and image code being finite 
essential initial factor codes.

\smallskip

\noindent {\small \rm (3)} Every right ideal morphism of $n\,A^*$, with 
domain code and image code being finite essential initial factor codes, is 
an extension of some element of $n{\cal RI}_A^{\sf fin}$.
\end{pro}
{\sc Proof.} (1) This follows from the examples of $F$, $F_1$, $F_2$ above; 
see also \cite[Lemma 2.27 and Fig.\ 2]{BinG}.

\smallskip

\noindent (2) Existence is trivial, since a joinless code is an initial 
factor code. To prove uniqueness we first prove two important properties:

\smallskip

\noindent [Finiteness Property]
 \ {\it A right ideal morphism $f$ whose ${\rm domC}(f)$ is finite and 
essential can only be extended in finitely many ways to a right ideal 
morphism of $nA^*$. }
 
Indeed, if ${\rm Dom}(f)$ is a finitely generated essential right ideal,
then $nA^* \minus {\rm Dom}(f)\,$ is finite. Hence, $f$ can be defined on
$\,nA^* \minus {\rm Dom}(f)\,$ in finitely many ways only.

\medskip

\noindent [Union Property] 
 \ {\it Let $f$ be a right ideal morphisms of $nA^*$ such that 
${\rm domC}(f)$ and ${\rm imC}(f)$ are finite essential initial factor codes.
Let $f_1$ and $f_2$ be any right ideal morphisms of $nA^*$ that extend $f$, 
such that ${\rm domC}(f_i)$ and ${\rm imC}(f_i)$ are finite essential initial 
factor codes (for $i = 1,2$). 
Then $f_1 \cup f_2$ is also a right ideal morphism of $nA^*$ with finite 
domain code and finite image code (that are essential initial factor codes), 
extending $f$. }

\smallskip

\noindent Part (2) of the Proposition follows immediately from the 
Finiteness and Union Properties. 

\smallskip

Let us prove the Union Property. Let $f_{12} = f_1 \cup f_2\,$ (the 
set-theoretic union). So 
${\rm Dom}(f_{12}) = {\rm Dom}(f_1) \cup {\rm Dom}(f_2)$, and this is an
essential right ideal (since ${\rm Dom}(f_1)$ and ${\rm Dom}(f_2)$ are 
essential right ideals). 
Then for $x \in {\rm Dom}(f_1) \minus {\rm Dom}(f_2)$ we have 
$f_{12}(x) = f_1(x)$; and for $x \in {\rm Dom}(f_2) \minus {\rm Dom}(f_1)$ 
we have $f_{12}(x) = f_2(x)$. This, and the next Claim, imply that 
$f_{12}$ is a function.

\medskip

\noindent {\sf Claim.} \ For all $x \in {\rm Dom}(f_1) \cap {\rm Dom}(f_2)$:
$ \ f_1(x) = f_2(x)$ \ ($\,= f_{12}(x)$).

\smallskip

\noindent Proof of the Claim: Let $x \in {\rm Dom}(f_1) \cap {\rm Dom}(f_2)$.
Since $f_1$ and $f_2$ are extensions of $f$, and ${\rm Dom}(f)$ is an 
essential right ideal, there exists $\ell \ge 0$ such that 
$\, x \cdot nA^{\ell} \subseteq {\rm Dom}(f)$.
Then for all $u \in nA^{\ell}$: $f_1(x) \, u = f_1(xu)$ $=$ $f(xu)$  $=$ 
$f_2(xu) = f_2(x) \, u$, hence $f_1(x)\,u = f_2(x)\,u$. Since $nA^*$ is a 
cancelative monoid it follows that $f_1(x) = f_2(x)$. 
 \ \ \ [End, Proof of the Claim.]

\smallskip

\noindent Now $f_{12}$ is uniquely defined on 
$\,{\rm Dom}(f_1) \minus {\rm Dom}(f_2)$, 
$\,{\rm Dom}(f_2) \minus {\rm Dom}(f_1)$, and 
$\,{\rm Dom}(f_1) \cap {\rm Dom}(f_2)$; hence it is a function on all of 
$\,{\rm Dom}(f_1) \cup {\rm Dom}(f_2)$.

It follows also that $f_{12}$ is a right ideal morphism. Indeed, on 
$\,{\rm Dom}(f_1)$ $\,=\,$
$({\rm Dom}(f_1) \minus {\rm Dom}(f_2))$ $\cup$
$({\rm Dom}(f_1) \cap {\rm Dom}(f_2))\,$ we have $f_{12} = f_1$. And on 
$\,{\rm Dom}(f_2)$ $\,=\,$
$({\rm Dom}(f_2) \minus {\rm Dom}(f_1))$ $\cup$
$({\rm Dom}(f_1) \cap {\rm Dom}(f_2))\,$ we have $f_{12} = f_2$. 

\smallskip

\noindent (3) Let $f: P\,nA^* \to Q\,nA^*$ be a right ideal morphism, where
$P, Q \subseteq nA^*$ are essential finite initial factor codes, and let 
$\,\ell = {\rm maxlen}(P)$.
We can restrict $f$ to the right ideal $n\,A^{\ge \ell}$, which is 
generated by the maximal joinless code $n\,A^{\ell}$. 
By Lemma \ref{LEMpropRI}, since $n\,A^{\ge \ell}$ is a maximal joinless 
code, $f(n\,A^{\ge \ell})$ is a maximal joinless code. 
 \ \ \ $\Box$

\bigskip

The {\em Brin-Thompson group} $n G_{k,1}$ is defined by
$ \ n \, {\cal RI}_A^{\sf fin}/\!\! \equiv_{\rm end}$.
 \ Equivalently (by \cite[Lemma 2.24]{BinG}): $\,n G_{k,1}$ is the action
monoid of $n \, {\cal RI}_A^{\sf fin}$ on $nA^{\omega}$. 
By the above Remark and Prop.\ \ref{LEMuniqext}, $n G_{k,1}$ is also the 
set of maximally extended right ideal morphisms of $nA^*$, where 
multiplication is composition followed by maximum extension. 
 \ When $k = 2$ we obtain $n V$.

\bigskip

The definitions of evaluation functions, evaluation problem, domain 
membership problem, and long versus short data inputs, can be generalized 
immediately to $2V$ over a finite generating set $\Gamma_{\!2}$.

\begin{defn} \label{DEFev2VGam} {\bf (evaluation problem of $2V$ over 
$\Gamma_{\!2}$).} 

\noindent The {\em evaluation problem} of $2V$ over a finite generating set 
$\Gamma_{\!2}$ is specified as follows. 

\noindent {\sc Input:} \ $(w, x, y) \in$ 
$\Gamma_{\!\!2}^* \x 2\{0,1\}^* \x 2\{0,1\}^*$.

\noindent {\sc Question:} \ $E_w(x) = y$ ?

\smallskip

\noindent Here, $E_w$ is the unique maximum extension of 
$w_n \circ \,\ldots\, \circ w_1(.)$ to a right ideal morphism, where
$w = w_n \ldots w_1$ with $w_n, \,\ldots\,, w_1 \in \Gamma_{\!2}$.
\end{defn}

\noindent For $x = (x_1,x_2) \in 2\,\{0,1\}^*$ we define
$\ell(x) = {\rm max}\{|x_1|, |x_2|\}$. Let 

\smallskip

 \ \ \  \ \ \ $\lambda_{\Gamma_2} = {\rm maxlen}(\Gamma_{\!2})$  $=$ 
$\max\{ \ell(z) : z \in \,\bigcup_{\gamma \in \Gamma_{\!2}}$
${\rm domC}(\gamma) \cup {\rm imC}(\gamma) \,\}$.  

\smallskip

\noindent By \cite[Cor.\ 3.3]{BinG} we have for all $f \in 2V$: 

\smallskip

 \ \ \  \ \ \  if $ \ \ell(x) \ge \lambda_{\Gamma_2}\,|w| \ $ then 
$\,w_n \circ \,\ldots\, \circ w_1(x)\,$ is defined.

\smallskip

\noindent Hence, similarly to Def.\ \ref{DEFlongshortdata}, for a given 
$w \in \Gamma_{\!\!2}^*\,$ we call a data input $x \in 2\,\{0,1\}^*$ a 
{\em long data input} for $w\,$  iff 
$\,w_n \circ \,\ldots\, \circ w_1(x)\,$ is defined.
A data input $x \in {\rm Dom}(E_w)$ that is not long is called a {\em short 
data input} for $w$.

\begin{defn} \label{DEFev2VGamLong}  
 \ The {\em evaluation problem} for $2V$ over $\Gamma_{\!2}$ for {\em long 
data inputs} is defined as follows.

\noindent {\sc Input:} \ $(w,x,y)$  $\in$
$\Gamma_{\!\!2}^* \x \{0,1\}^* \x \{0,1\}^*$, where $w = w_n\,\ldots\,w_1$,
with $w_n, \,\ldots\,, w_1 \in \Gamma_{\!1}$.

\noindent {\sc Question:} \ $w_n \circ \,\ldots\, \circ w_1(x) = y$ ?
\end{defn}

\section{Results}

We now look at the complexity of the evaluation problem for $V$ and $2V$, 
and we compare the evaluation problem with the word problem.
The evaluation decision problem for $V$ and $2V$, as well as the word 
problem, depend on $\Gamma_{\!1}$, respectively $\Gamma_{\!2}$.
However, the complexity changes only slightly with the generating set,
provided that $\Gamma_{\!1}$ and $\Gamma_{\!2}$ are finite.

\begin{pro} \label{VoverGamma} {\bf (evaluation problem of {\boldmath $V$ 
over $\Gamma_{\!1}$ for long data inputs}).}

\noindent Let $\Gamma_{\!1}$ be a finite generating set of the Thompson group 
$V$.  The evaluation problem of $V$ over $\Gamma_{\!1}$ for {\em long} data 
inputs is in $\,{\sf DCF} \cap {\sf DCF}^{\rm rev}$.  More precisely (assuming 
$\,\Gamma_{\!1} \cap \{0,1\}$  $=$  $\varnothing$ and 
${\mathbb 1} \in \Gamma_{\!1}$), the language

\smallskip

 \ \ \   \ \ \  $L_{_{\!V}} \ = \ \{x^{\rm rev} w\, y \,: \ $
$x, y \in \{0,1\}^*, \ w \in \Gamma_{\!\!1}^+,$
$ \ w_n \circ \,\ldots\, \circ w_1(x) = y\,\}$

\smallskip

\noindent is deterministic context-free, and its reverse is also 
deterministic context-free.
\end{pro}
{\sc Proof.} Let $w = w_n \ldots w_1$ with $w_n, \,\ldots\,, w_1$  $\in$
$\Gamma_{\!1}$. The set $\{0,1\}^* \, \Gamma_{\!\!1}^+ \, \{0,1\}^*$ is 
finite-state, so we can check easily whether the input has the correct 
format. 

We construct a deterministic push-down automaton (dpda) that on input
$x^{\rm rev} w y$ proceeds as follows. First, the dpda reads $x^{\rm rev}$
and pushes it onto the stack; the left end of $x$ (i.e., the right end of
$x^{\rm rev}$) is now at the top of the stack. 
Next, the dpda reads $w_1$, pops the prefix $p_1$ of $x$ that belongs to 
${\rm domC}(w_1)$, applies the function $w_1$ to $p_1$, and pushes 
$w_1(p_2)$ onto the stack. This pop-push cycle is repeated with $w_2$, 
$\,\ldots\,$, $w_n$.
Any one of these steps could be undefined; $x$ is a long data input iff
all steps are defined. So in this process, the dpda automatically checks 
whether $x$ is a long data input for $w$.
When $w$ has been entirely read, the stack content (if defined) is
$w(x)$; here, the left end of $w(x)$ is at the top of the stack.
Next, the dpda compares the stack content with $y$.
It accepts iff the stack is empty at the moment the reading of $y$ ends.

It is well known that the class {\sf DCF} is not closed under reversal. 
But $\,L_{_{\!V}}^{\rm \,rev} \,=\, \{y^{\rm rev} w^{\rm rev} x :\,$
$x, y \in \{0,1\}^*, \ w \in \Gamma_{\!\!1}^+,$
$ \ w_n \circ \,\ldots\, \circ w_1(x) = y\,\}$
 \ is accepted in the same way as $L_{_{\!V}}$, except that now the dpda
checks the relation $w_1^{-1} \circ \,\ldots\, \circ w_n^{-1}(y) = x$ 
(which is equivalent to $w_n \circ \,\ldots\, \circ w_1(x) = y$).
 \ \ \  $\Box$

\medskip

In Prop.\ \ref{VoverGamma} the syntax of the input is important; i.e., the 
input is assumed to be of the form $x^{\rm rev} w\, y\,$ (not $x w y$, nor 
$y w x$, etc.). For example, the set
$\{x w x : x \in \{0,1\}^*, \,w \in \Gamma_{\!\!1}^+\}$ is not in {\sf CF}, 
but $\{x^{\rm rev} w x : x \in \{0,1\}^*, \,w \in \Gamma_{\!\!1}^+\}$ is in
{\sf DCF}.

\begin{thm} \label{VoverCiruit} {\bf (evaluation problem of 
{\boldmath $V$ over $\Gamma_{\!1} \cup \tau$, and $2V$ over $\Gamma_{\!2}$},
for long data inputs).} 

\smallskip

The evaluation problem of the Thompson group $V$ over a circuit-like 
generating set $\Gamma_{\!1} \cup \tau$, for {\em long} data inputs, is 
{\sf P}-complete with respect to log-space many-one reduction.

The evaluation problem of the Brin-Thompson group $2V$ over a finite 
generating set $\Gamma_{\!2}$, for {\em long} data inputs, is 
{\sf P}-complete, with respect to log-space many-one reduction.
\end{thm}
Theorem \ref{VoverCiruit} will be proved in Section 6.

\bigskip

\bigskip

\noindent {\Large \bf Relation between the evaluation problem and the word
problem}

\medskip

We will compare the evaluation problems of $V$ and $2V$ with the word
problem.
The word problem of $V$ over a finite generating set is in {\sf coCF}
(Lehnert and Schweitzer \cite{LS}); $\,{\sf DCF} \cap {\sf DCF}^{\rm rev}$ 
and {\sf DCF} are strict subclasses of ${\sf CF} \cap {\sf coCF}$ and 
{\sf coCF}.  
On the other hand, the word problem of $V$ over a circuit-like generating 
set $\Gamma_{\!1} \cup \tau$, and the word problem of $2V$ and $nG_{k,1}$ 
over a finite generating set, are {\sf coNP}-complete \cite{BinG, BinGk1}.

\smallskip

Although the evaluation problem and the word problem of $V$ and $2V$ look 
similar to the circuit value problem, respectively the circuit equivalence 
problem, the following proposition shows that there is also a fundamental 
difference. This is caused by the existence of short (versus long) data
inputs.

\begin{pro} \label{VoverGammatoWP} {\bf (reduction of the word problem
to the evaluation problem).}

\smallskip

\noindent {\small \rm (1)} Let $\Gamma_{\!1}$ be a finite generating set of 
the {\em Thompson group} $V$. 
The evaluation problem of $V$ over $\Gamma_{\!1}$, or over the circuit-like 
generating set $\Gamma_{\!1} \cup \tau$, for the data input and output 
$\varepsilon$, is equivalent to the word problem of $V$ over $\Gamma_{\!1}$,
respectively $\Gamma_{\!1} \cup \tau$.
More precisely, for any $w$ in $\Gamma_{\!\!1}^*$ or in 
$(\Gamma_{\!1} \cup \tau)^*$,

\smallskip

 \ \ \   \ \ \  $E_w(\varepsilon) = \varepsilon$ \ \ \ iff 
 \ \ \ $w = {\mathbb 1}\,$ in $V$.

\smallskip

\noindent Hence the evaluation problem (in general) of $V$ over 
$\Gamma_{\!1}$ is in {\sf coCF}; and the evaluation problem of $V$ over 
$\Gamma_{\!1} \cup \tau$ is {\sf coNP}-complete.

\medskip

\noindent {\small \rm (2)} The evaluation problem of the 
Brin-Thompson group $2V$ over a finite generating set $\Gamma_{\!2}$, for the
data input and output $(\varepsilon, \varepsilon)$, is equivalent to the
word problem of $2V$. 
More precisely, for any $w \in \Gamma_{\!\!2}^*$, 

\smallskip

 \ \ \   \ \ \  $E_w((\varepsilon, \varepsilon))$  $=$ 
$(\varepsilon, \varepsilon)$ \ \ \ iff
 \ \ \ $w = {\mathbb 1}\,$ in $2V$.

\smallskip

\noindent Hence the evaluation problem of $2V$ over $\Gamma_{\!2}$
is {\sf coNP}-complete.

\medskip

\noindent {\small \rm (3)} 
The word problem of $V$ over a finite generating set $\Gamma_{\!1}$, or over
$\Gamma_{\!1} \cup \tau$, can be reduced to the evaluation problem of $V$ 
over $\Gamma_{\!1}$, respectively  $\Gamma_{\!1} \cup \tau$, for data inputs 
of length $N$, for any fixed $\,N \le O(\log |w|)$; the reduction is a 
polynomial-time conjunctive reduction of polynomial arity. 

Similarly, the word problem of $2V$ over a finite generating set 
$\Gamma_{\!2}$ can be reduced to the evaluation problem of $2V$ over 
$\Gamma_{\!2}$ data inputs $(x_1, x_2) \in 2\,\{0,1\}^*$ of length 
$\,\max\{|x_1|,\, |x_2|\} = N$, for any fixed $\,N \le O(\log |w|)$;
the reduction is a polynomial-time conjunctive reduction of polynomial arity.
\end{pro}
{\sc Proof.} 
(1), (2): The equivalences are straightforward, since $\{\varepsilon\}$ is a 
maximal prefix code in $\{0,1\}^*$, and $\{(\varepsilon, \varepsilon)\}$ is 
a maximal joinless code in $2\,\{0,1\}^*$.

\noindent (3) We observe that $w = {\mathbb 1}$ in $V$ iff $w(x) = x$ for all 
$x \in$ $\{0,1\}^N$, for any $N \in {\mathbb N}$. Indeed, $\{0,1\}^N$ is a
maximal prefix code. If we pick $N \le O(\log |w|)$ then the cardinality of 
$\{0,1\}^N$ is bounded by a polynomial in $|w|$, and all $x \in \{0,1\}^N$ 
can be found in polynomial time.

A similar reasoning applies to $2V$, since $2\,\{0,1\}^N$ is a maximal
joinless code.
 \ \ \  $\Box$

\begin{thm} \label{redEVtoWP} {\bf (reduction of the evaluation problem
to the word problem).}

\smallskip

\noindent {\small \rm (1)} The evaluation problem of $V$ over 
the circuit-like generating set $\Gamma_{\!1} \cup \tau$ is in {\sf coNP},
and can be reduced to the word problem of $V$ over $\Gamma_{\!1} \cup \tau$ 
by a many-one $\log$-space reduction.
 
Similarly, the evaluation problem of $2V$ over a finite generating set 
$\Gamma_{\!2}$ is in {\sf coNP}, and can be reduced to the word problem of 
$2V$ over $\Gamma_{\!2}$ by a many-one $\log$-space reduction. 

\medskip

\noindent {\small \rm (2)} The evaluation problem of $V$ over a finite 
generating set $\Gamma_{\!1}$ (and more generally, over a generating set
$\Gamma_{\!1} \cup \Delta$), reduces to the word problem of $V$ over 
$\Gamma_{\!1}$ (respectively $\Gamma_{\!1} \cup \Delta$) by an 
eight-fold conjunctive $\log$-space reduction.  Here $\Delta$ is any infinite 
subset of $V$, coded by a $\log$-space set of strings.
\end{thm}
{\sc Proof.} (1) The evaluation problem is described by the  
following {\sf coNP}-formula:

\smallskip

 \ \ \ $E_w(x) = y$ \ \ iff \ \ $(\forall z \in \{0,1\}^*)[\,|z|$  $=$ 
$c_{_{\Gamma_1, w}} |w| - |x|$ $\,\Rightarrow\,$  $E_w(xz) = yz\,]$,

\smallskip

\noindent where $c_{_{\Gamma_1, w}}$  $=$ 
$\max\{c_{_{\Gamma_1}},\, {\rm maxindex}_{\tau}(w)\}$,
as seen before.

The length of the string $z$ in the $\forall$-quantifier is 
polynomially bounded; indeed, $|z| \le c_{_{\Gamma_1, w}} \, |w| - |x|$,
and $c_{_{\Gamma_1, w}}$ is linearly bounded in terms of the size of $w$.

The predicate $w(xz) = yz$ can be checked in polynomial time. Indeed, by 
the condition $|xz| = c_{_{\Gamma_1, w}}\,|w|$, the data input $xz$ is 
{\em long}.
Hence $E_w(xz) = w_n \circ \,\ldots\, \circ w_1(xz)$, so $E_w(xz)$ is simply
computed (in at most quadratic time) by applying the generators $w_i$ in 
sequence. 
Hence the universal formula above is a {\sf coNP}-formula, so the problem 
is in {\sf coNP}. 

The word problem of $V$ over $\Gamma_{\!1} \cup \tau$ is {\sf coNP}-complete,
hence the evaluation problem (being in {\sf coNP}) reduces to the word 
problem.

\smallskip

For $2V$ over $\Gamma_{\!2}$, the same reasoning works; here we use the fact 
that $\ell(x) = \max\{x_1, x_2\} \le \lambda_{\Gamma_2} \, |w|$ implies that
$x$ is a long data input. 

\smallskip

\noindent (2) This proof is given in Section 7. 
 \ \ \  \ \ \ $\Box$

\bigskip

\bigskip

\noindent {\bf Open problems:}

\smallskip

\noindent $\bullet$ \ Is the word problem of $M_{2,1}$ over a finite 
generating set {\sf P}-complete? \ (It is in {\sf P} by 
\cite{BiThompsMonV3}.) 

\smallskip

\noindent $\bullet$ \ Does there exist a {\em finitely presented} group, 
or monoid, whose word problem over a finite generating set is 
{\sf P}-complete? 
 \ \ (Compare with other {\sf P}-complete problems for groups in 
\cite[A.8.6 - A.8.17]{GHR}.) 

\bigskip

\bigskip

\newpage

\noindent {\bf Overview of the reductions between problems}

\smallskip

\noindent The following graph shows the reductions between the problems 
considered in this paper. An arrow $A \to B$ indicates that $A$ reduces to 
$B$ by many-one log-space reduction or conjunctive log-space reduction.
Mutual reduction is indicated by $A \leftrightarrow B$.

\bigskip

\noindent Abbreviations:

\smallskip

wp.\ \ --- \ the word problem for a given group and generating set

ev.\ \ --- \ the evaluation problem for a given group, generating set, and 
data input

$\Gamma_{\!1}$ \ --- \ any finite generating set of $V$

$\Gamma_{\!2}$ \ --- \ any finite generating set of $2V$

$\Delta$ \ --- \ any subset of $V$ that is encoded in binary by a 
log-space language

\bigskip

\bigskip

\bigskip

\bigskip

\noindent
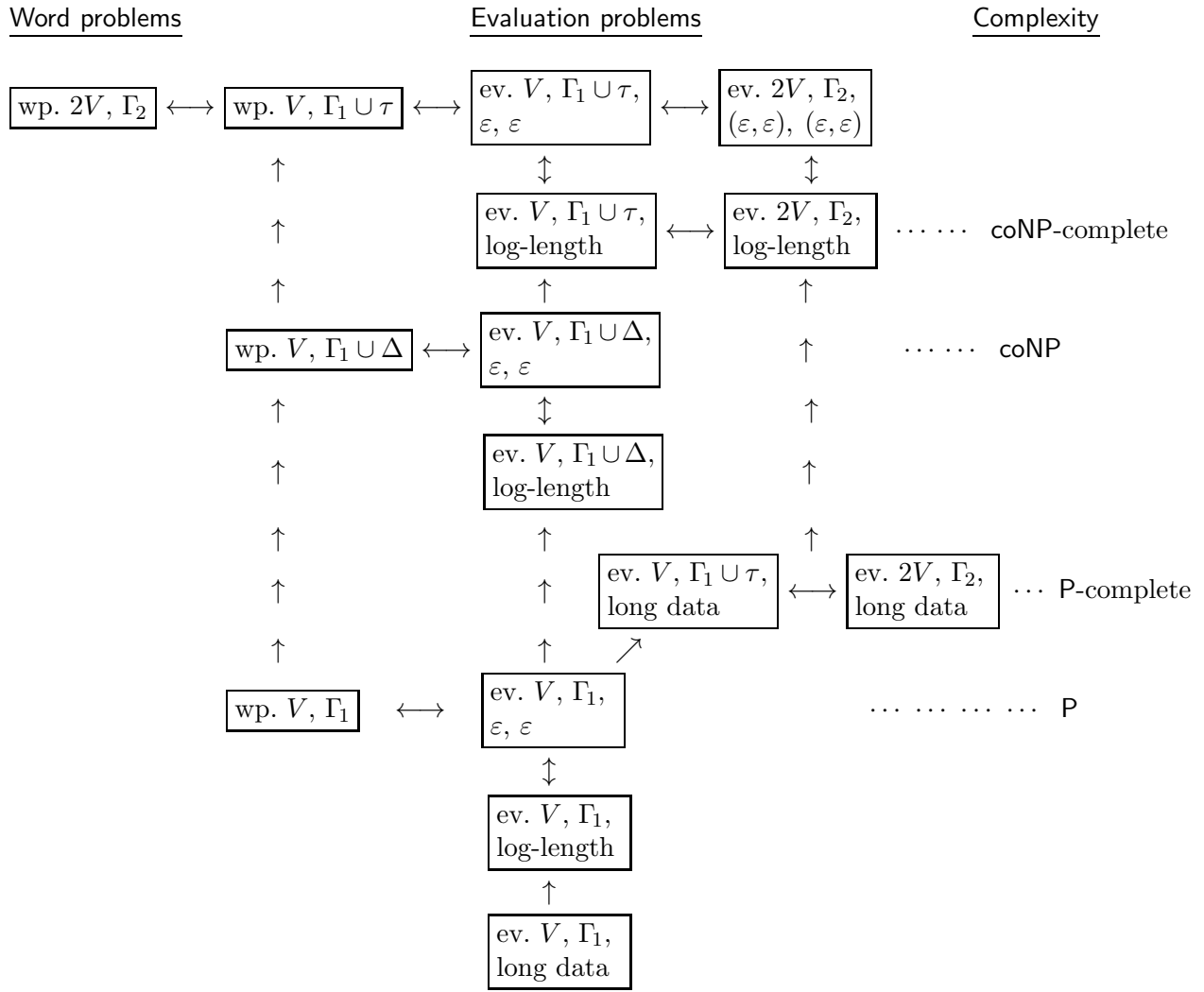
\begin{minipage}{\textwidth}

\noindent \underline{\sf Word problems}  \hspace{1.5in}  
\underline{\sf Evaluation problems}  \hspace{1.4in}
\underline{\sf Complexity}

\bigskip

\fbox{wp.\ $2V$, $\Gamma_{\!2}$}
$\longleftrightarrow$ 
\fbox{wp.\ $V$, $\Gamma_{\!1} \cup \tau$} 
$\longleftrightarrow$
\fbox{\parbox{0.9in}{ev.\ $V$, $\Gamma_{\!1} \cup \tau$,\\   
      $\e$, $\e$}}
$\longleftrightarrow$
\fbox{\parbox{0.75in}{ev.\ $2V$, $\Gamma_{\!2}$,\\  
                     $(\e,\e), \ (\e,\e)$}}

\smallskip

\hspace{1.4in} $\uparrow$ \hspace{1.3in}  $\updownarrow$
\hspace{1.3in} $\updownarrow$

\smallskip

\hspace{1.4in} $\uparrow$  \hspace{0.96in}
\fbox{\parbox{0.9in}{ev.\ $V$, $\Gamma_{\!1} \cup \tau$,\\   
       log-length}}
$\longleftrightarrow$
\fbox{\parbox{0.75in}{ev.\ $2V$, $\Gamma_{\!2}$,\\
       log-length}}
 \ $\cdots \, \cdots$ \ {\sf coNP}-complete

\smallskip

\hspace{1.4in} $\uparrow$ \hspace{1.3in}  $\uparrow$
\hspace{1.27in} $\uparrow$

\smallskip
 
\hspace{1.15in} 
\fbox{wp.\ $V$, $\Gamma_{\!1} \cup \Delta$}
$\longleftrightarrow$
\fbox{\parbox{0.9in}{ev.\ $V$, $\Gamma_{\!1} \cup \Delta$,\\   
      $\e$, $\e$}}
\hspace{0.67in} $\uparrow$
\hspace{0.4in} $\cdots \, \cdots$ \ {\sf coNP}

\smallskip

\hspace{1.4in} $\uparrow$ \hspace{1.3in} $\updownarrow$ \hspace{1.3in} $\uparrow$

\smallskip

\hspace{1.4in} $\uparrow$ \hspace{1.0in} 
\fbox{\parbox{0.9in}{ev.\ $V$, $\Gamma_{\!1} \cup \Delta$,\\
       log-length}} \hspace{0.67in} $\uparrow$

\medskip

\hspace{1.4in} $\uparrow$ \hspace{1.3in} $\uparrow$ 
 \hspace{1.3in} $\uparrow$

\hspace{1.4in} $\uparrow$ \hspace{1.3in} $\uparrow$ \hspace{0.16in} 
\fbox{\parbox{0.9in}{ev.\ $V$, $\Gamma_{\!1} \cup \tau$,\\  
     long data}} 
$\longleftrightarrow$
\fbox{\parbox{0.75in}{ev.\ $2V$, $\Gamma_{\!2}$,\\
     long data}} $\,\cdots\,$ {\sf P}-complete

\smallskip

\hspace{1.4in} $\uparrow$ \hspace{1.3in} $\uparrow$ 
\hspace{0.25in} $\nearrow$

\smallskip

\hspace{1.15in}
\fbox{wp.\ $V$, $\Gamma_{\!1}$}
\hspace{0.1in} $\longleftrightarrow$ \hspace{0.1in}
\fbox{\parbox{0.7in}{ev.\ $V$, $\Gamma_{\!1}$,\\
        $\e$, $\e$}}
\hspace{1.2in}  \ $\cdots \ \cdots \ \cdots \ \cdots$ \ {\sf P}

\smallskip

\hspace{2.9in} $\updownarrow$ 

\smallskip

\hspace{2.6in}
\fbox{\parbox{0.7in}{ev.\ $V$, $\Gamma_{\!1}$,\\
         log-length}}

\smallskip

\hspace{2.9in}  $\uparrow$

\smallskip

\hspace{2.6in}  
\fbox{\parbox{0.7in}{ev.\ $V$, $\Gamma_{\!1}$,\\
       long data}} 

\end{minipage}


\newpage

\section{Complementary prefix codes, partial fixators, commutation 
tests}

This Section introduces tools to be used in Section 6 for proving
Theorem \ref{VoverCiruit}, and in Section 7 for proving Theorem 
\ref{redEVtoWP}(2).
Let $A$ be a finite alphabet of cardinal $\,|A| = k \ge 2$.

\subsection{Complementary prefix codes}

\begin{defn} \label{DEFcomplprefc} {\rm \cite[Def.\ 5.2]{BiCoNP}.}
 \ Two prefix codes $P,P' \subseteq A^*$ are {\em complementary prefix codes} 
 \ iff \  $\,P \cup P'\,$ is a maximal prefix code in
$A^*$, and $\,PA^* \cap P' A^* = \varnothing$.
\end{defn}
This definition is equivalent to the following: $\,P \cup P'\,$ is a maximal
prefix code, and $\,P \cap P' = \varnothing$.

If $P$ is a maximal prefix code then $\varnothing$ is the unique 
complementary prefix code of $P$; except for this case, the complementary
prefix code of $P$ is never unique. E.g., for every $u \in P'$, 
 \ $(P' \minus \{u\}) \cup uA\,$ is also a complementary prefix code of $P$.

\bigskip

\noindent {\bf Notation:} For an alphabet $A$ and a set $S \subseteq A^*$,
let 

\smallskip

 \ \ \  \ \ \ ${\sf pref}(S)$  $\,=\,$ 
$\{x \in A^* : (\exists s \in S)[\,x \le_{\rm pref} s\,]\}$; 

\smallskip

\noindent i.e., ${\sf pref}(S)$ is the set of prefixes of elements of $S$.
And 

\smallskip

 \ \ \  \ \ \  ${\sf Spref}(S)$  $=$ 
$\{x \in A^* : (\exists s \in S)[x < _{\rm pref} s]\}$; 

\smallskip

\noindent i.e., ${\sf Spref}(S)$ is the set of strict prefixes of elements 
of $S$.

\begin{pro} \label{PROPexistcomplem} {\bf (existence and construction).}
 \ For every finite prefix code $P \subseteq A^*$ with $|A| \ge 2$,
there exists a finite complementary prefix code $P' \subseteq A^*$. If 
$P$ is not $\varnothing$ and not a maximal prefix code, then 
${\rm maxlen}(P') = {\rm maxlen}(P)$.

If $P$ is given by a list of strings, then a complementary prefix code $P'$ 
can be computed from $P$ in $\log$-space.
\end{pro}
{\sc Proof.} This is a special case of
\cite[Lemma 2.27 and Cor.\ 2.30]{BinMk1}. If $P$ is maximal then 
$P' = \varnothing$, and if $P = \varnothing$ then we can choose 
$P' = \{\e\}$. 
If $P$ is not $\varnothing$ and not a maximal prefix code, then a
complementary prefix code of $P$ is given by the formula

\medskip

 \ \ \  \ \ \  $P' \,=\,$   $\{xa \,: \ $
$x \in {\sf Spref}(P), \ a \in A, \ xa \not\in {\sf pref}(P)\}$.

\medskip

\noindent Claim 1: \ $P'$ is a prefix code. 

\smallskip

\noindent Let us consider any $xa, yb \in P'$, where 
$x, y \in {\sf Spref}(P)$ and $a,b \in A$. 

Case 1:
If $x$ and $y$ are prefixes of a same $p = p_1 \ldots p_k \in P$ then 
$x \,\|_{\rm pref}\, y$; we can assume $x \le_{\rm pref} y$.
If $x = y$ and $a \ne b$ then $xa$ and $yb = xb$ differ in position $|x|+1$,
so they are not prefix-comparable.
If $x <_{\rm pref} y$, then $x = p_1 \ldots p_i$, and  
$y = p_1 \ldots p_i p_{i+1} \ldots p_j$, for some $1 \le i < j < k$.
Since $xa \not\in {\sf pref}(P)$,  $a \ne p_{i+1}$ so $xa$ differs from $yb$
in position $i+1$; hence $xa$ and $yb$ are not prefix-comparable.

Case 2: If $x, y$ are not in Case 1, and if
$x <_{\rm pref} p$, $y <_{\rm pref} q$, for $p, q \in P$, 
then $x$ and $y$ are not prefix-comparable (otherwise $x, y$ are in Case 1).
Then $x$ and $y$ differ in some position $i \le \min\{|x|, |y|\}$, so $xa$ 
and $yb$ also differ in that position; hence $xa$ and $yb$ are not 
prefix-comparable.  \ [End, Claim 1.]

\smallskip

\noindent Claim 2: \ $P \cap P' = \varnothing$.

The formula says that $xa \not\in {\sf pref}(P)$, hence $xa \not\in P$,
hence $P \cap P' = \varnothing$. 

\smallskip

\noindent Claim 3: \ $P \cup P'$ is maximal.

Let $z \in A^*$ be any string of length $|z| = {\rm maxlen}(P)$ such that
$z \not\in P A^*$.  Then exists $s \in {\sf Spref}(P)$ such the $s$ is a 
prefix of $z\,$ (where $s$ could be $\e$). 
If $s = z_1 \ldots z_i$ is the longest prefix of $z$ that belongs to 
${\sf pref}(P)$ then $s z_{i+1}$ $\not\in$ ${\sf pref}(P)$ with
$i = |s|\,$ (otherwise we would have $z \in PA^*$). 
Therefore, $s z_{i+1}$ $\in$ $P'$.  Thus for every $z \in A^*$ with 
$|z| = {\rm maxlen}(P)$: either $z \in P A^*$, or $z \in P' A^*$. 
 \ [End, Claim 3.]

\smallskip

\noindent Claims 1, 2, 3 imply that $P'$ is a complement of $P$. 

\smallskip

Thanks to this formula, $P'$ can be computed from $P$ in log-space, as 
follows. First, by comparing all strings in $P$ one by one, it can be 
checked that $P$ is a prefix code. Non-emptiness is trivial to check. 
(Non)maximality can be checked by the Kraft (in)equality.
After that, for every $p \in P$, every strict prefix $s$ can be found, hence
each string $sa$ (for $a \in A$) is found, after which one can check whether 
$sa$ is not a prefix of some element in $P$.
The output of this process is a list of all the elements of $P'$ in an 
arbitrary order, with possible repetitions.
Since sorting can be done in log-space, $P'$ can be listed in increasing 
dictionary order, without repetitions. 
Recall that the composite of log-space functions is in log-space 
\cite[Lemma 13.3]{HU}.
 \ \ \  \ \ \ $\Box$

\bigskip

\noindent A special case is of Prop.\ \ref{PROPexistcomplem} is the 
following.

\begin{pro} \label{PROPexistcomplemSingl}
 \ For every $u = u_1 \ldots u_{\ell} \in A^+$ with 
$u_1, \,\ldots\,, u_{\ell} \in A$ and $\ell = |u| > 0$,
the following is a {\em complementary prefix code} of $\{u\}$:

\medskip

 \ \ \   \ \ \  $\ov{u} \ = \ $
$\bigcup_{j=0}^{\ell-1} \,u_1 \ldots u_j \, (A \minus \{u_{j+1}\})$;

\medskip

\noindent here $\,u_1 \ldots u_j = \e\,$ when $j = 0$.
 \ The set $\ov{u}$ is a prefix code with the following properties:

\smallskip

\noindent $\bullet$ \ \ ${\rm maxlen}(\ov{u}) = |u|$,

\smallskip

\noindent $\bullet$ \ \ $|\ov{u}| = |u| \, (|A| - 1)$, 

\smallskip

\noindent $\bullet$ \ \ $\ov{u}$ can be computed (as a list of strings) 
from $u$ in time $O(|u|^2)$, and in $\log$-space. 
 \ \ \   \ \ \  $\Box$
\end{pro}
If $u = \varepsilon$ then $\{u\}$ is a maximal prefix code, so
the complementary prefix code is $\varnothing$.

\begin{lem} \label{LEMxTOy} {\bf (transitivity of $V$).} 

\noindent Suppose $V$ is defined as the set of maximally extended right ideal
morphisms in ${\cal RI}_A^{\sf fin}$.

For any $u, v \in \{0,1\}^+$ there exists $\psi \in V$ such that 
$\psi(u) = v$ and such that the table size of 
$\psi$ is $ \ |{\rm domC}(\psi)| \,=\, 1 + \max\{|u|, |v|\}$. 
\end{lem}
{\sc Proof.} Let $\ov{u}$ and $\ov{v}$ be complementary prefix codes of
$\{u\}$, respectively $\{v\}$, as in Prop.\ \ref{PROPexistcomplemSingl}.
So, the cardinalities satisfy $|\ov{u}| = |u|$ and $|\ov{v}| = |v|$. 

If $|u| > |v|$, we replace the prefix code $\ov{v}$ by 
$Q_v = (\ov{v} \minus \{x\}) \cup \{x0, x1\}$, for some $x \in \ov{v}$;
this preserves the fact that $Q_v$ is a complementary prefix code of 
$\ov{u}$, but increases its cardinality by 1. We repeat this until 
$|Q_v| = |\ov{u}| = |u|$. The case where $|v| > |u|$ is handled in a similar 
way. In any case, we obtain complementary prefix codes $Q_u$ and $Q_v$ of 
$\{u\}$, respectively $\{v\}$, such that $|Q_u| = |Q_v| = \max\{|u|, |v|\}$.

Now, $\psi$ is defined by $\psi(u) = v$ and $\psi$ maps $Q_u$ bijectively 
onto $Q_v$, arbitrarily. The table size of $\psi$ is $1 + \max\{|u|, |v|\}$.
 \ \ \   \ \ \ $\Box$

\subsection{Partial fixators}

We generalize the well known concept of fixator to partial injections.

\begin{defn} \label{DEFpFix}
{\bf (partial fixator} {\rm \cite[Def.\ 4.13]{BinG}}{\bf ).}

\smallskip

\noindent A function $g$ {\em partially fixes} a set $S$ \ iff
 \ $g(x) = x\,$ for every $\, x \in $
$S \,\cap\, {\rm Dom}(g) \,\cap\, {\rm Im}(g)$.
The {\em partial fixator} of $S$ in $V$ is

\medskip

 \ \ \ \ \  ${\rm pFix}_V(S)$  $ \ = \ $  $\{g \in V : \,$
$(\forall x \in S \cap {\rm Dom}(g) \cap {\rm Im}(g))[\,g(x) = x\,]\,\}$.
\end{defn}
This is also called partial pointwise stabilization.

We will only use partial fixators for sets $S$ that are right ideals of 
$\{0,1\}^*$. If $S = P\,\{0,1\}^*$ is a right ideal, where $P$ is a prefix 
code, then ${\rm pFix}_V(S)$ is a group \cite[Lemma 4.1]{BiCoNP}.
We abbreviate $\,{\rm pFix}_V(P\,\{0,1\}^*)\,$ by $\,{\rm pFix}_V(P)$.
In particular, for $z \in \{0,1\}^*$ we abbreviate
${\rm pFix}_V(z\,\{0,1\}^*)$ by ${\rm pFix}_V(z)$.

One easily proves that if $P = \varnothing$ then $\,{\rm pFix}_V(P) = V$,
and if $P$ is a maximal prefix code then
$\,{\rm pFix}_V(P) = \{{\sf id}\}\,$ (the one-element subgroup of $V$).

\begin{lem} \label{Fixes}
 \ Below, $\subseteq$ or $=$ refer to set inclusion or equality of subsets 
of $V$.  \ For all $u, v \in \{0,1\}^*$:  
 
\smallskip

 \ \ \ ${\sf pFix}_V(u) \,\subseteq\, {\sf pFix}_V(v)$
 \ \ iff \ \ $v\, \{0,1\}^* \subseteq u\, \{0,1\}^*$ \ (i.e., 
$\,u \le_{\rm pref} v$).  

\smallskip

\noindent Hence, 

\smallskip

 \ \ \ ${\sf pFix}_V(u) \,=\, {\sf pFix}_V(v)$ \ \ iff \ \ $u = v$.
\end{lem}
{\sc Proof.} The equality relations immediately follow from the inclusions.
Let us prove the inclusions.

\smallskip

\noindent $[\Leftarrow]$ If $v = uz \in u\,\{0,1\}^*$, and if $f$ fixes 
$u$ then $f$ fixes $uz = v\,$ (since $f$ is a right ideal and $f(u)$ is 
defined).

\smallskip

\noindent $[\Rightarrow]$ We prove the contrapositive: Suppose 
$u \not\le_{\rm pref} v$.

\smallskip

\noindent Case 1: $u$ and $v$ are not prefix-comparable.

Then $\{u,v\}$ is a finite prefix code, so there exists a finite 
complementary prefix code $Q$.

If $Q = \varnothing$ then $\{u,v\} = \{0,1\}$, since this is the only 
two-element maximal prefix code in $\{0,1\}^*$. 
Suppose $u = 0$ (the case $u = 1$ is similar). 
Then $f$ with table $\{(0,0),\, (10,11)\, (11, 10)\}$ satisfies $f \in$ 
${\sf pFix}_V(u) \minus {\sf pFix}_V(v)$. \ So, 
${\sf pFix}_V(u) \not\subseteq {\sf pFix}_V(v)$.

If $Q \ne \varnothing$, let $q \in Q$. Then $f$ with table 
$\,\{(u,u),\, (v,q),\, (q,v)\} \,\cup\, \{(z,z) : z \in Q \minus \{q\}\}\,$
satisfies $f \in$  ${\sf pFix}_V(u) \minus {\sf pFix}_V(v)$.  \ So, 
${\sf pFix}_V(u) \not\subseteq {\sf pFix}_V(v)$.

\smallskip

\noindent Case 2: $u$ and $v$ are prefix-comparable.

Then $v <_{\rm pref} u\,$ (since $u \not\le_{\rm pref} v$). So, $u = va$ for
some $a \in \{0,1\}^+$. We can assume $a = 0b$ for some $b \in \{0,1\}^*$
(the case $a = 1b$ is similar); so, $u = v0b$. 
Then $\{u, v10, v11\}$ is a prefix code; let $Q$ be a complementary prefix 
code.  Then $f$ with table 
$\{(u,u),\, (v10, v11),\, (v11,v10)\} \,\cup\, \{(z,z) : z \in Q\}$
satisfies $f \in {\sf pFix}_V(u) \minus {\sf pFix}_V(v)$.  \ So, 
${\sf pFix}_V(u) \not\subseteq {\sf pFix}_V(v)$.
 \ \ \ $\Box$

\bigskip

\noindent The next theorem is a generalization of \cite[Lemma 4.20]{BinG}.

\begin{thm} \label{FixP} {\bf (generators of ${\sf pFix}_V(P)$).}
 \ For any non-empty, non-maximal, finite prefix code $P \subseteq \{0,1\}^*$
we have:  

\smallskip

\noindent {\small \rm (1)} \ \ ${\sf pFix}_V(P)\,$ is isomorphic to $V$  
 \ (hence, ${\sf pFix}_V(P)\,$ is finitely generated). 

\smallskip

\noindent {\small \rm (2)} \ \ Let $\Gamma_{\!1}$ be a finite generating set
of $V$; $\,\Gamma_{\!1}$ will be kept fixed.  
Then, from any $P$ (given as a list of strings), a finite generating set for
$\,{\sf pFix}_V(P)$ ($\,\subseteq V$) can be computed in $\log$-space, with 
each generator of ${\sf pFix}_V(P)$ expressed as a word over 
$\,\Gamma_{\!1}$.
\end{thm}
{\sc Proof.} (1) Let $Q = \{q_j : 1 \le j \le k\}$ be a finite complementary 
prefix code of $P$. Let $B = \{b_1, \ldots, b_k\}$ be an alphabet of cardinal
$\,|B| = k = |Q|$, such that $B \cap \{0,1\} = \varnothing$, and let us 
consider the Higman-Thompson group $G_{2,k}$. Every element of $G_{2,k}$ is 
given by a bijection between two maximal prefix codes in $B\,\{0,1\}^*$. 
For $\psi \in G_{2,k}$, the domain and image codes are of the form 
$\,{\rm domC}(\psi)$ $\,=\,$ $\bigcup_{i=1}^k b_i C_i$, and
$\,{\rm imC}(\psi)$ $\,=\,$ $\bigcup_{j=1}^k b_j D_j$, 
where each $C_i$ and $D_j$ is a finite maximal prefix code in 
$\{0,1\}^*$ (for $i, j = 1, \ldots, k$); and $\psi$ is given by a bijection 
from ${\rm domC}(\psi)$ onto ${\rm imC}(\psi)$; so,
$|{\rm domC}(\psi)| = |{\rm imC}(\psi)|$. 

To construct an isomorphism from $G_{2,k}$ onto ${\sf pFix}_V(P)$ we choose 
a bijection $f_{B,Q}: B \to Q$. For simplicity, we assume that $Q$ is 
indexed in increasing dictionary order, and we pick $f_{B,Q}(b_i)$ to 
be the $i$th element of $Q$ in the dictionary order of $\{0,1\}^*$; thus
$\,f_{B,Q}(b_i) = q_i$. An isomorphism from $G_{2,k}$ onto ${\sf pFix}_V(P)$ 
is now defined by

\medskip

 \ \ \  $\Pi: \ \psi \in G_{2,k}$  $ \ \longmapsto \ $ 
$\Pi(\psi)$  $\,=\,$ 
${\sf id}_P \,\cup\, \pi(\psi)$ \  $\,\in\,$  ${\sf pFix}_V(P)$,

\medskip

\noindent where $\,{\rm domC}(\pi(\psi)) \,=\, \bigcup_{i=1}^k q_i C_i$, 
$\,{\rm imC}(\pi(\psi))$ $\,=\,$ $\bigcup_{j=1}^k q_j D_j$, and

\medskip

 \ \ \  $\pi(\psi)\big(q_i c_i\big) = q_j d_j$ \ iff 
     \ $\psi(b_i c_i) = b_j d_j$.

\medskip

\noindent Then $\pi(\psi)$ is well-defined; indeed, since $Q$ is a prefix 
code, $q_i c_i$ uniquely determines $i$, $b_i$ and $c_i$; hence $q_j d_j$ 
is uniquely determined by $\psi$ and $q_i c_i$. 

So, ${\rm domC}(\Pi(\psi)) = P \ \cup \  \bigcup_{i=1}^k q_i C_i$, and
$\,{\rm imC}(\Pi(\psi))$ $\,=\,$ $P \ \cup \ \bigcup_{j=1}^k q_j D_j$; 
these are finite maximal prefix codes. Each of 
$ \ \bigcup_{i=1}^k q_i C_i \ $ and $ \ \bigcup_{j=1}^k q_j D_j \ $
is a complementary prefix code of $P$.
It is straightforward to see that $\Pi$ is invertible and that it is a 
homomorphism.  So, $G_{2,k}$ is isomorphic to ${\sf pFix}_V(P)$.

Finally, by Higman's \cite[Cor.\ 2, p.\ 12]{Hig74}, all the free 
J\'onson-Tarski algebras ${\cal V}_{2,k}$ are isomorphic to ${\cal V}_{2,1}$. And since $G_{2,k}$ is the group of automorphisms of ${\cal V}_{2,k}$, all 
the groups $G_{2,k}$ are isomorphic to $G_{2,1}$ for all $k \ge 1$. 
Hence, $G_{2,1}$ is isomorphic to ${\sf pFix}_V(P)$.
 
\medskip

\noindent (2) Recall that $\Gamma_{\!1}$ is fixed and $P$ is a variable 
input; for this input $P$, our goal is to compute generators of 
${\sf pFix}_V(P)$, expressed as words over $\Gamma_{\!1}$. For this, we 
first construct an isomorphism $\,\Pi_{_P}: V \to {\sf pFix}_V(P)\,$ as in
part {\small \rm (1)}, but we describe this isomorphism more explicitly. 
(The construction in {\small \rm (1)} is thus redundant, but {\small \rm (1)}
makes {\small \rm (2)} easier to understand.)

Let $B \subseteq \{0,1\}^*$ be any maximal prefix code such that $|B| = |Q|$;
therefore, now $B$ is a set of bitstrings, as opposed to the construction in
{\small \rm (1)}, where $B$ was a new alphabet. Also, $B$ is a maximal prefix
code, whereas $Q$ is non-maximal.
Let $f_{B,Q}: B \to Q\,$ be a bijection. For simplicity, we pick $f_{B,Q}$ 
so as to preserve the dictionary order of $\{0,1\}^*$; then 
$\,f_{B,Q}(b_i) = q_i\,$ (for $i = 1,\ldots,k$), assuming that $B$ and $Q$ 
are both indexed in increasing dictionary order.

For any $\varphi \in G_{2,1}$ we restrict $\varphi$ so that 
$\,{\rm domC}(\varphi) \,\cup\, {\rm imC}(\varphi)$  $\,\subseteq\,$
$B \{0,1\}^*$; we still call this restriction $\varphi$. 
Then $\,{\rm domC}(\varphi) = \bigcup_{i=1}^k b_i C_i\,$ and 
$\,{\rm imC}(\varphi) = \bigcup_{j=1}^k b_j D_j\,$
for some finite maximal prefix codes $C_i, D_j \subseteq \{0,1\}^*\,$ 
(for $i, j = 1, \ldots, k$).
We define an isomorphism from $G_{2,1}$ onto ${\sf pFix}_V(P)$ by

\medskip

 \ \ \  \ \ \  $\Pi_{_P} : \ \varphi \in G_{2,1}$  $ \ \longmapsto \ $
$\Pi_{_P}(\varphi)$  $\,=\,$  ${\sf id}_P \,\cup\, \pi(\varphi)$  
 \ \  $\,\in\,$  ${\sf pFix}_V(P)$,

\medskip

\noindent where $ \ {\rm domC}(\pi(\varphi))$  $\,=\,$ 
$\bigcup_{i=1}^k q_i C_i$,
$ \ {\rm imC}(\pi\varphi))$ $\,=\,$ $\bigcup_{j=1}^k q_j D_j$; \ and
$\pi(.)$ is such that

\medskip

 \ \ \  \ \ \  $\pi(\varphi)\big(q_i c_i\big) = q_j d_j$ \ \ iff
  \ \ $\varphi(b_i c_i) = b_j d_j$.

\medskip

\noindent We can also write

\medskip

 \ \ \  \ \ \  $\Pi_{_P}(\varphi)(.)$  $\,=\,$  ${\sf id}_P$  $ \ \cup \ $
$f_{B,Q} \circ \varphi \circ f_{B,Q}^{\, -1}(.)$,

\medskip

\noindent where $f_{B,Q}: B \to Q\,$ is the bijection that preserves the
dictionary order, seen above.
Then $\pi(\varphi)$ is well-defined. Indeed, since $Q$ is a prefix
code, $q_i c_i$ uniquely determines $q_i$ and $c_i$; and $q_i$ determines 
$i$ and $b_i$. And $b_i c_i$ determines $b_j d_j$ ($\,= \varphi(b_i c_i)$),
from which $b_j$ and $d_j$, hence $q_j = \pi(b_j)$, hence $q_j d_j$, 
are determined.
As in part {\small \rm (1)}, we see that $\Pi_{_P}$ is an isomorphism from
$G_{2,1}$ onto ${\sf pFix}_V(P)$.

\medskip

Let us now find generators for ${\sf pFix}_V(P)$.
Since $\Pi_{_P}(.)$ is an isomorphism from $G_{2,1}$ onto ${\sf pFix}_V(P)$, 
it follows immediately that $\Pi_{_P}(\Gamma_{\!1})$ is a finite generating 
set of ${\sf pFix}_V(P)$, that is in one-to-one correspondence with
$\Gamma_{\!1}$. 
In the computation of $\Pi_{_P}(\Gamma_{\!1})$ from $P$, only $P$ is variable 
($\Gamma_{\!1}$ is fixed); so the complexity of the calculation of 
$\Pi_{_P}(\gamma)$, for each $\gamma \in \Gamma_{\!1}$, is measured as a 
function of $\,\sum_{p \in P} |p|$.
Since $\Gamma_{\!1}$ is fixed and finite, we can look at $\Pi_{_P}(\gamma)$ 
for one $\gamma \in \Gamma_{\!1}$ at a time.  

The composite of two log-space computable functions is log-space computable 
\cite[Lemma 13.3]{HU}. 
Therefore, in order to compute a word for $\Pi_{_P}(\gamma) \in$ 
$\Gamma_{\!\!1}^*$ from $P$ in log-space, it suffices to show that each one 
of the following steps is computable in log-space:

\smallskip

\noindent $\circ \ $  from $P$, find a complementary prefix code $Q$;

\noindent $\circ \ $  from $Q$, find a maximal prefix code $B$  $\subseteq$
$\{0,1\}^*$ with $\,|B| = |Q|$;

\noindent $\circ \ $  from $B$ and $Q$, find 
$\,{\rm domC}(\Pi_{_P}(\gamma))\,$ and $\,{\rm imC}(\Pi_{_P}(\gamma))\,$ 
(by also using ${\rm domC}(\gamma)$ and ${\rm imC}(\gamma)$);

\noindent $\circ \ $  from ${\rm domC}(\Pi_{_P}(\gamma))$ and
${\rm imC}(\Pi_{_P}(\gamma))$, find the table of $\Pi_{_P}(\gamma)\,$ (using 
$\gamma$ and the formula for $\Pi_{_P}(.)$);

\noindent $\circ \ $  from the table of $\Pi_{_P}(\gamma)$, find a word for 
$\Pi_{_P}(\gamma) \in \Gamma_{\!\!1}^*$ such that this word evaluates to
$\Pi_{_P}(\gamma)$; this 

\!word has length 
$ \ \le \ c\,\big(|\Pi_{_P}(\gamma)| \ \log |\Pi_{_P}(\gamma)| \big)$ 
 \ for some constant $c > 0$.

\medskip

\noindent $\bullet \ $ Finding $Q$: Prop.\ \ref{PROPexistcomplem} gives a 
formula for the complementary prefix code,

\smallskip

 \ \ \  $Q$ $=$  $\{xa :\, $
$x \in {\sf Spref}(P), \ a \in \{0,1\}, \ xa \not\in {\sf pref}(P)\}$,

\smallskip

\noindent from which $Q$ can be computed in log space (as a list of strings,
sorted in increasing dictionary order).  

\smallskip

\noindent $\bullet \ $  Finding $B$: We  choose 
$\,B = \{0^{k-1}\} \,\cup\, 0^{\le k-2} 1$, where $k = |Q|$.  The set $B$ 
(in the form of a list of strings) can be found from $Q$ in log-space. 
%

\smallskip

\noindent $\bullet \ $  Finding the table for $\Pi_{_P}(\gamma)$ from the
table of $\gamma$ is done in two steps: 
First we find the table of $\gamma_B$ (the restriction of $\gamma$ to $B$); 
second, from this we find the table of $\Pi_{_P}(\gamma)$.

Finding the table of $\gamma_B$: We look at ever string $x \in$ 
${\rm domC}(\gamma)$. If $x$ has a prefix in $B$, then $(x,\,\gamma(x))$ is 
already part of the table of $\gamma_B$, so we simply copy it to the output.
If $x$ is a strict prefix of some string $b \in B$, i.e., $b = xu$ for some
$u \in \{0,1\}^+$, then we include $(b, \,\gamma(x) \, u)$ into the table of
$\gamma_B$.
Recall also that by \cite[Lemma 3.3]{BiThomps},
$\,{\rm domC}(\gamma_B) \subseteq B \,\cup\, {\rm domC}(\gamma)$, so 
$\gamma_B$ has table-size
$\,|{\rm domC}(\gamma_B)| \le |B| + |{\rm domC}(\gamma)|$.
 \ Since ${\rm domC}(\gamma)$ and $B$ are given as lists of strings, the 
above procedure can be carried out in log-space: Only a bounded set of 
positions in the input needs to be kept track of. 

\smallskip

Finding the table of $\Pi_{_P}(\gamma)$ from the table of $\gamma_B$ and $B$:
The input is $\,B$, $\gamma$, $\gamma_B$, and $\pi: B \to Q\,$ (all given as 
lists of strings or pairs of strings). 
By definition, the table of $\Pi_{_P}(\gamma)$ consists of 
$ \ \{(p,p) : \in P\}$ $\,\cup\,$
$\,\{(q_ic_i,\, q_j d_j) : (b_ic_i,\, b_j d_j) \in \gamma_B\}$.  

Outputting $\,\{(p,p) : \in P\}\,$ is easy.

For the second part, we look at each $(b_ic_i,\, b_j d_j)$ in the table of
$\gamma_B\,$ (for $i = 1,\ldots,k$). 
From $b_ic_i$ and the prefix code $B$, we find $b_i$; by using the table of
$\pi$ we find and output $q_i$; from $b_ic_i$ and the prefix code $B$ we 
also find and output $c_i$; so we have produced ``$(q_i c_i,$''.
In a similar, from $b_j d_j$ we find and output ``$q_j d_j)$''.
All this can be done in log-space, since only a bounded set of positions in 
the input need to kept track of.

\smallskip

\noindent $\bullet \ $  Finding a word over $\Gamma_{\!1}$ for 
$\Pi_{_P}(\gamma)$: This can be done in log-space by using the construction 
in \cite[Proof of Theorem 3.8(1)]{BiThomps}. This word has length 
$\,\le$ $O\big(|\Pi_{_P}(\gamma)| \, \log |\Pi_{_P}(\gamma)| \big)$ 
$=$ $O(k \, \log k)$.

\smallskip

Notation: For any $g \in V$, $\|g\|$ denotes the table size of $g$
(i.e., the number of pairs in the table). The input size for $g$, used for 
measuring complexity, is $ \ \sum_{x \in {\rm domC}(g)} |x|$  $\,+\,$ 
$\sum_{y \in {\rm imC}(g)} |y|$.
 
\smallskip

Let us examine the the construction in 
\cite[Proof of Theorem 3.8(1)]{BiThomps}, where for any $g \in V$, given by 
its table, a word for $g$ over $\Gamma_{\!1}$, of length 
$O(\|g\| \, \log \|g\|)$, is computed in space 
$\,O\big(\log (\sum_{x \in {\rm domC}(g)} |x|$  $+$ 
$\sum_{y \in {\rm imC}(g)} |y|) \big)$. This done in several steps, as
follows.

\smallskip

\noindent $*$ \ A maximal prefix code $S_n$ of size $n = \|g\| \,$ is 
chosen, such that $S_n \subseteq \{0,1\}^{k-1} \cup \{0,1\}^k$, where 
$k = \lceil \log n \rceil$. See \cite[Prop.\ 3.9 and following]{BiThomps}.  

Let $F$ denote the well-known Thompson group of dictionary order
preserving elements of $V$.
Now $g$ is factored as $g = \beta_g\,\pi_g\,\alpha_g$, where
$\beta_g, \alpha_g \in F$, and ${\rm domC}(\pi_g) = {\rm imC}(\pi_g) = S_n$.
The three factors are given by their tables; see \cite[Prop.\ 3.9]{BiThomps}. 
This is easily done in log-space.

\smallskip

\noindent $*$ \ The elements $\beta_g, \alpha_g$ are factored over the 
two well-known generators of $F$; see \cite[Prop.\ 3.10]{BiThomps}.
This can be implemented in log-space.

\smallskip

\noindent $*$ \ The element $\pi_g$ is factored into $\,\le 3n$ 
string transpositions of the form 
$\,(0^{\lceil \log n \rceil} \ | \ w)$, where $0^{\lceil \log n \rceil}$  
$\in$  $S_n$, and $w \in S_n \minus \{0^k\}$.  This factorization is easily
carried out in log-space, based on the table of $\pi_g$.

Here we define string transpositions $(0^k\,|\,w)$ as follows. For any 
$w \in \{0,1\}^+$, let $j \in {\mathbb N}$ be the unique number such that 
$0^j 1$ is a prefix of $w$, so $w = 0^j1 v$.  Now we let

\smallskip

\hspace{0.1in} $(0^k|w) \,=\, \{(0^k,w),\, (w,0^k)\}$  

\hspace{0.8in} $\,\cup\,$  
$\{(x,x) :\, x = 0^i 1, \ 0 \le i < k, \ i \ne j\}$  

\hspace{0.8in} $\,\cup\,$ 
$\{(x,x) :\, x = 0^j 1 p a, \ p \in \{0,1\}^*,\,$
$p <_{\rm pref} v, \ a \in \{0,1\}, \ p a \not\le_{\rm pref} v \}$.

\smallskip

\noindent $*$ \ Every transposition $(0^k|w)$ is factored over 
$\Gamma_{\!1}$ as a word of length $\,O(\log n)$.
See \cite[Lemma 3.11]{BiThomps}. The proof of this Lemma leads to 
four computational steps: 

Cases (1.1) and (1.2) can be handled in one step, using log-space.
 
Cases (1.3) and (1.4) can be handled in one step, using log-space.

Case (2.1) can be handled in one step, using log-space; this leads to
case (2.2).

Case (2.2) can be handled in one step, using log-space; this leads to 
cases (2.3).

Cases (2.3) and (2.4) can be handled jointly, using log-space. 

\smallskip

\noindent The result is a word over $\Gamma_{\!1}$ whose value in $G_{2,1}$ 
is $\,\Pi_{_P}(\gamma)$. 
 \ \ \  \ \ \ $\Box$

\subsection{Commutation test}

Let $A$ be a finite alphabet of cardinal $|A| = k \ge 2$, and let $G_{k,1}$ 
be the Higman-Thompson group. 
The following ``commutation test'' for membership in 
${\sf pFix}_{G_{k,1}}(P)$ was introduced in \cite[Cor.\ 9.15]{BiCoNP}.

\begin{thm} \label{CommTestP} {\bf (commutation test for membership in
${\sf pFix}_{G_{k,1}}(P)$).} 

\noindent Let $P, Q \subseteq A^*$ be any finite, non-empty, 
{\em complementary} prefix codes; let $\Gamma_{\!Q}$ be any generating set 
of ${\sf pFix}_{G_{k,1}}(Q)$.
For every $g \in G_{k,1}$ the following are equivalent: 

\smallskip

\noindent {\small \rm (0)} \hspace{0.2in} 
    $g \,\in\, {\sf pFix}_{G_{k,1}}(P)$,

\smallskip 

\noindent {\small \rm (1)} \hspace{0.2in}
$g$ commutes with every element of $\,{\sf pFix}_{G_{k,1}}(Q)$,

\smallskip

\noindent {\small \rm (2)} \hspace{0.2in}
$g$ commutes with every generator in $\,\Gamma_{\!Q}$.
\end{thm}
{\sc Proof.} It is straightforward to see that (1) and (2) are equivalent.
For the proof that (0) is equivalent to (1), see \cite[Cor.\ 9.15]{BiCoNP}.
 \ \ \ $\Box$

\medskip

\noindent The commutation test can be effectively used when $\Gamma_Q$ 
is finite.

\section{The evaluation problem of $V$ over $\Gamma_{\!1} \cup \tau$, and for
$2V$ over $\Gamma_{\!2}$, for long data inputs}

We will prove that these two problems are {\sf P}-complete (Theorem 
\ref{VoverCiruit}). Let $A$ be any finite alphabet with $|A| = k \ge 2$.
As before, let $\Gamma_{\!1}$ and $\Gamma_{\!2}$ be a finite generating set
of $V$, respectively $2V$.

\begin{lem} \label{LEMevVcirc2VinP} {\bf (the problems are in {\sf P}).}

\noindent The evaluation problem of $V$ over $\Gamma_{\!1} \cup \tau$, and 
the evaluation  problem for $2V$ over $\Gamma_{\!2}$, for {\em long} data 
inputs, belong to {\sf P}.
\end{lem}
{\sc Proof.} For $w \in (\Gamma_{\!1} \cup \tau)^*$ and $x, y \in \{0,1\}^*$,
we apply $w_1$ to $x$, then $w_2$ to $w_1(x)$, $w_3$ to $w_2(x)$, etc.; this
is straightforward. 
The length of each $w_i \circ \,\ldots\, \circ w_1(x)\,$ (for 
$i = 1, \ldots, n$) is at most $\,c_{_{\Gamma_1, w_i\,\ldots\,w_1}} \cdot i$. 
Here, for any string $v \in (\Gamma_{\!1} \cup \tau)^*$:  
$\,c_{_{\Gamma_1, v}} = $
$\max \{ c_{_{\Gamma_1}},\, {\rm maxindex}_{\tau}(v) \}$, where
$c_{_{\Gamma_1}} = {\rm maxlen}(\Gamma_1)$, i.e., the length of the longest
bit-string in
$ \ \bigcup_{\gamma \in \Gamma_{\!1}} {\rm domC}(\gamma)$ $\cup$ 
${\rm imC}(\gamma)$;  and ${\rm maxindex}_{\tau}(v) = $
$\max \{i \in {\mathbb N}_{>0} : \tau_{i-1, i}$ occurs in $v \}\,$
(i.e., the largest subscript of any element of $\tau$ that occurs in $v$).

A very similar reasoning applies to $2V$ over $\Gamma_{\!2}$.
 \ \ \ $\Box$

\begin{lem} \label{LEMevVcircTo2V} {\bf (reduction from $V$ over 
$\Gamma_{\!1} \cup \tau$ to $2V$ over $\Gamma_{\!2}$).}

\noindent The evaluation problem of $V$ over $\Gamma_{\!1} \cup \tau$ can be
reduced to the evaluation problem for $2V$ over $\Gamma_{\!2}$,  by a 
many-one finite-state reduction. 
(This holds for all data inputs, long or short.)
\end{lem}
{\sc Proof.} Similarly to \cite[Section 4.6]{BinG}, we reduce the evaluation 
problem $(w,x,y)$ of $V$ over $\Gamma_{\!1} \cup \tau$, to the evaluation 
problem $\,(W,\, (x,\e),\, (y,\e))\,$ for $2V$ over $\Gamma_{\!2}$, where 
$W$ is obtained from $w$ as follows:

\smallskip

\noindent $\bullet$ \ every $\gamma \in \Gamma_{\!1}$ is replaced by 
$\,\gamma \x {\mathbb 1}$, defined by
$\,\gamma \x {\mathbb 1}(u,v) = (\gamma(u), v)$, for all $(u,v) \in$ 
$2\,\{0,1\}^*$, 

\smallskip

\noindent $\bullet$ \ every $\tau_{i,i+1} \in \tau$ occurring in $w$ is 
replaced by
$ \ \sigma^{i-1} \cdot (\tau_{1,2} \x {\mathbb 1}) \cdot \sigma^{-i+1}(.)$.

\smallskip

\noindent Then we have $\,E_w(x) = y\,$ iff $\,E_W((x,\e)) = (y,\e)$.

Moreover, $W$ can be computed from $w$ in log-space.
Outputting $\gamma \x {\mathbb 1}$ when $\gamma$ is read can be done by a
finite automaton. Outputting 
$\sigma^{i-1} \cdot (\tau_{1,2} \x {\mathbb 1}) \cdot \sigma^{-i+1}$
when $\tau_{i,i+1}$ is read can be done in log-space, assuming that
$\tau_{i,i+1}$ is encoded in binary by $a b^{i+1} a$. (See the remark about 
encoding in Section 1.)

So far we have used the following finite subset $\,\Gamma_{\!0}$  $\,=\,$ 
$\,\{\gamma \x {\mathbb 1} : \gamma \in \Gamma_{\!1}\}$
$\cup$ $\{\sigma, \ \tau_{1,2} \x {\mathbb 1}\,\}$ 
of $2V$ in order to generate all of $\tau$.
Next, every element of the finite set $\Gamma_{\!0}$ can be replaced by a
string over $\Gamma_{\!2}$; this can be done by a finite automaton.
  \ \ \ $\Box$

\bigskip

Finally, to prove Theorem \ref{VoverCiruit} we need to show that the 
evaluation problem of $V$ over $\Gamma_{\!1} \cup \tau$ with long data 
inputs is {\sf P}-hard. We do this by reducing the circuit value problem to 
this problem, using a many-one log-space reduction.
The difficulty is, of course, that we need to simulate arbitrary 
transformations by bijective functions.
For this we adapt the methods in \cite{BinG} (see also \cite{BiCoNP}), where 
the circuit equivalence problem was reduced to the word problem of $V$ over 
$\Gamma_{\!1} \cup \tau$.

\begin{defn} \label{simulation} {\rm \cite[Def.\ 4.10]{BinG}} 
{\bf (simulation):} 
 \ Let $\, f: \{0,1\}^m \to \{0,1\}^n \,$ be a total function. An element
$\Phi_f \in V$ {\em simulates} $f$ \ \ iff \ \ for all $\, x \in \{0,1\}^m$:
 \ \ \ $\Phi_f(0 \, x) \ = \ 0 \ f(x) \ x$.

When $\Phi_f$ is represented by a {\em word} $w_f \in$
$(\Gamma_{\!1} \cup \tau)^*$ \ we also say that $w_f$ simulates $f$.
\end{defn}

We want to define a {\em size} for every $w \in (\Gamma_{\!1} \cup \tau)^*$,
where $\Gamma_{\!1}$ is a finite generating set of $V$. First, 
${\sf size}(\gamma) = 1$ for every $\gamma \in \Gamma_{\!1}$, and 
${\sf size}(\tau_{i,i+1}) = i+1$ for every $\tau_{i,i+1} \in \tau$.
Finally, for $w = w_n \,\ldots\, w_1 \in (\Gamma_{\!1} \cup \tau)^*$ with
$w_n, \,\ldots\,, w_1 \in \Gamma_{\!1} \cup \tau$, the {\em size} of $w$ is 
$\, {\sf size}(w) \,= \ \sum_{i = 1}^n {\sf size}(w_i)$.

\smallskip

For a directed acyclic graph, the {\em depth} of a vertex $v$, denoted by 
${\sf depth}(v)$, is the length of a longest directed path ending in $v$. 
For any $d \in {\mathbb N}$, and any directed acyclic graph $C$, the set of
vertices that have depth $d$ is called the $d$th {\em layer}.  
An acyclic circuit (and more generally, a directed acyclic graph) is called 
{\em strictly layered} iff for every vertex $v$: all the parents of $v$ have 
the same depth ${\sf depth}(v) - 1$

For {\sf P}-completeness of the circuit value problem we only need to use 
acyclic circuits that are strictly layered. Indeed, the proof of 
{\sf P}-completeness of the circuit-value problem in 
\cite[Thm.\ 8.1]{Papadim} gives a $\log$-space many-one reduction of the 
acceptance problem of any one-tape polynomial-time Turing machine $M$ with a 
data input $x$ to a circuit-value problem $(C,x)$, where $C$ is strictly 
layered; see also \cite[Theorem 6.2.5]{GHR}, where strictly layered circuits 
are called synchronous. Therefore, \cite[Theorem 4.12]{BinG} can be 
replaced by the following, which gives a lower complexity.

\begin{thm} \label{reduction} {\bf (existence of a simulation).}
 \ There is an injective function  \ $C  \mapsto w_C$ \ from the set of
strictly layered acyclic boolean circuits to the set of words over 
$\Gamma_{\!1} \cup \tau$, with the following properties:

\smallskip

\noindent {\rm (1)} \ \ $w_C$ simulates the input-output function $f_C$
of $C$;

\smallskip

\noindent {\rm (2)} \ \ the size of $w_C$ satisfies
 \ \ $\|w_C\| \,<\, c \ |C|^3$ \ \ (for some constant $c > 0$);

\smallskip

\noindent {\rm (3)} \ \ $w_C$ is computable from $C$ in polynomial time,
in terms of $|C|$.
\end{thm}
{\sc Proof.} The proof appears in \cite[Theorem 4.12]{BinG}, where we can 
skip the last step, since $C$ is already strictly layered.
 \ \ \ $\Box$

\bigskip

Now the circuit value problem for strictly layered circuits (which is 
{\sf P}-complete) can be reduced to the evaluation problem of $V$ over 
$\Gamma_{\!1} \cup \tau$ as follows: 

For $(C,x,y)$, where $C$ is a strictly layered acyclic circuit, and $x, y$
$\in$ $\{0,1\}^*$,  consider $(w_C,\, 0x,\, 0yx)$ $\in$  
$(\Gamma_{\!1} \cup \tau)^* \x \{0,1\}^* \x \{0,1\}^*$. By Theorem 
\ref{reduction} we have:

\smallskip

\hspace{0.5in} $C(x) = y$ \ \ iff \ \ $w_C(0x) = 0yx$.

\smallskip

\noindent The function $(C,x,y) \mapsto (w_C,\, 0x,\, 0yx)$ is a many-one 
polynomial-time reduction, by Theorem \ref{reduction}. The proof of 
\cite[Theorem 4.12]{BinG} also shows that $w_C(0x)$ is defined when the 
generators in $w_C$ are applied to $0x$; i.e., $0x$ is a {\em long} data input 
for $w_C$.

\bigskip

Finally, let us prove that the reduction $\,(C,x,y)$  $\mapsto$ 
$(w_C,\, 0x,\, 0yx)\,$ is actually a many-one {\em log-space} reduction (and
not just polynomial-time). 
For this, we review the proof of \cite[Theorem 4.12]{BinG}, except that we
skip the parts where the circuit is transformed into a strictly layered
circuit.
For $1 < j$, let 

 \ \ \  \ \ \  $\sigma_{1,j}$  $\,=\,$
$\tau_{j-1,j} \tau_{j-2,j-1} \,\ldots\, \tau_{2,3} \tau_{1,2}(.)$;

\smallskip

\noindent in other words, $\sigma_{1,j}$ is the cyclic shift
 \    $a_j a_{j-1} \,\ldots\, a_{2} a_1$ 
$\,\longmapsto\,$  $a_1 a_j a_{j-1} \,\ldots\, a_{2}$.

Let $C$ be a strictly layered circuit with $L$ layers, with data input 
$\,x_1 \,\ldots\, x_m$  $\in$  $\{0,1\}^*$, and data output 
$\,y_1 \,\ldots\, y_n$  $\in$  $\{0,1\}^*$. 
Let $C_{\ell}$ be the circuit consisting of layer $\ell$ (for 
$1 \le \ell \le L$). 
The output of $C_{\ell}$ is the bitstring $Y^{\ell}$; its input is 
$Y^{\ell-1}$. 
The circuit $C_{\ell}$ is simulated by the word $w_{C_{\ell}}$ over 
$\Gamma_{\!1} \cup \tau$; the function represented by $w_{C_{\ell}}$ is
$\,\Phi_{C_{\ell}} :$  $0 Y^{\ell -1} \longmapsto 0 Y^{\ell} Y^{\ell -1}$.

The word $w_{C_{\ell}}$ is constructed by simulating each gate in 
$C_{\ell}$. E.g., if $w_{C_{\ell}}$ contains an {\sc and} gate in position
$i+1$ (counting the gates from left to right), with input variables 
$x'_{i+1}, x'_{i+2}$ (in $Y^{\ell-1}$) and output variable $y'_{j+1}$ (in
$Y^{\ell}$), then

\smallskip
 
 \ \ \  \ \ \  $w_{C_{\ell}}$  $\,=\,$  
$v_{(C_{\ell}, > i+1)}$  $\sigma_{1,j+2}$  
$\tau_{3,j+i+4} \tau_{2,j+i+3}$  $\varphi_{\wedge}$ 
$\tau_{3,j+i+4} \tau_{2,j+i+3}$  $\varphi_{\rm 0f}$   
$v_{(C_{\ell}, < i+1)}\,$.

\smallskip

\noindent 
Here $v_{(C_{\ell},\,>i+1)} \in (\Gamma_{\!1} \cup \tau)^*$ simulates 
the gates to the right of gate number $i+1$, and $v_{(C_{\ell},\,<i+1)}$ 
simulates the gates to the left of gate $i+1$; $\varphi_{\wedge}$ simulates 
the {\sc and}-gate, and $\varphi_{\rm 0f}$ simulates the {\sc fork}-gate (see
\cite[following Lemma 4.11]{BinG}). For the other types of gates (namely 
{\sc or}, {\sc not}, and {\sc fork}), the representations by a word over 
$\Gamma_{\!1} \cup \tau$ is similar (see \cite[proof of Thm.\ 4.12]{BinG}).

The entire circuit $C$ is simulated by the word  

\smallskip

 \ \ \  \ \ \  $w_C$  $=$ 
$\pi_2$  $(w_{C_{L-1}} \ \ldots \ w_{C_1})^{-1}$  $\pi_1$
$w_{C_L} w_{C_{L-1}} \ \ldots \ w_{C_1}$,

\smallskip

\noindent representing the function

\smallskip

 \ \ \  \ \ \  $\Phi_C: 0 x_1 \,\ldots\, x_m$  $\longmapsto$ 
$0 y_1 \,\ldots\, y_n x_1 \,\ldots\, x_m$.

\smallskip

\noindent Here, $\pi_1 = (\sigma_{1,|Z|})^n$, where 
$Z = 0 y_1 \,\ldots\, y_n Y^{L-1} \,\ldots\, Y^2 Y^1 x_1 \,\ldots\, x_m$, so
$|Z| = 1 +n+ m + \sum_{\ell = 1}^{L-1} |Y^{\ell}|$ (i.e., 1 plus the input
length, plus the sum of the output lengths of all the layers).
And $\pi_2 = (\sigma_{1,n+m})^m$.

\smallskip

To compute $w_C$ from $C$, a Turing machine just needs to keep track of 
positions inside the circuit $C$ (namely the layer $\ell$, and the position 
within the layer). This can be done in log space.

This is straightforward to see for $\pi_2$, since from $C$ one can directly 
find the number $m$ of input wires, and the number $n$ of output wires. 

To output $\,(w_{C_{L-1}} \ \ldots \ w_{C_1})^{-1}$ ($\, =$
$w_{C_1}^{-1} \ \ldots \ w_{C_{L-1}}^{-1}$), the Turing machine 
can compute each $w_{C_{\ell}}$ (for $\ell = L-1, \ldots, 2, 1$). The 
subcircuit $C_{\ell}$ is directly obtained from $C$; and $w_{C_{\ell}}$ is 
the concatenation of the representations of all the gates occurring in 
$C_{\ell}$, from left to right. 
E.g., for an {\sc and} gate at position $i+1$ within $C_{\ell}$ with output
wire at position $j+1$, the representation of the gate (seen above) is 
$\,\sigma_{1,j+2}$
$\tau_{3,j+i+4} \tau_{2,j+i+3}$  $\varphi_{\wedge}$
$\tau_{3,j+i+4} \tau_{2,j+i+3}$  $\varphi_{\rm 0f}$;
this can be computed in log-space, based on $i$ and $j$.

To output $\pi_1$, the Turing machine needs $|Z|$, which is a number obtained
from the size of $C$, hence in log-space. Then $\pi_1 = (\sigma_{1,|Z|})^n$
can be directly found. 

Finally, $w_{C_L} w_{C_{L-1}} \ \ldots \ w_{C_1}$ is found in a similar way
as $\,(w_{C_{L-1}} \ \ldots \ w_{C_1})^{-1}$. 

From $\,w_C$  $=$
$\pi_2$  $(w_{C_{L-1}} \ \ldots \ w_{C_1})^{-1}$  $\pi_1$
$w_{C_L} w_{C_{L-1}} \ \ldots \ w_{C_1}$,
a log-space Turing machine can obtain a word over $\Gamma_{\!1} \cup \tau$, 
by eliminating inverses; this mainly involves  reordering the word (since the
transpositions in $\tau$ are their own inverses). Recall that we assume that 
$\Gamma_{\!\!1}^{-1} = \Gamma_{\!\!1}$. Recall that the composite of 
log-space computable functions is log-space computable \cite[Lemma 13.3]{HU}.
 \ \ \ $\Box$

\section{Reduction of evaluation problems to word problems}

\subsection{Reduction to the monoid word problem}

The circuit value problem is easily reduced to the equivalence problem of 
circuits. A monoid version $M_{2,1}$ of the Thompson group $G_{2,1}$ was
defined in \cite{BiThompsMon, BiThompsMonV3}. 
The evaluation problem for the Thompson monoid $M_{2,1}$ (over a 
finite generating set $\Gamma_M$ or a circuit-like generating set 
$\Gamma_M \cup \tau$) can be reduced to the word problem of $M_{2,1}$ 
($\Gamma_M$, respectively $\Gamma_M \cup \tau$), as follows. 

For any strings $u, v \in \{0,1\}^*$ let $\,[v \leftarrow u](.)\,$ denote
the element of $M_{2,1}$ with table $\{(u,v)\}$. In particular, 
$[u \leftarrow u] = {\sf id}_{u \{0,1\}^*}$, i.e.\ the identity function 
restricted to $u\,\{0,1\}^*$. And $[\varepsilon \leftarrow u](.)$ is the
{\em pop} $u$ operation (erasing a prefix $u$, and undefined on inputs that 
do not have $u$ as a prefix); and $[v \leftarrow \varepsilon](.)$ is the 
{\em push} $v$ operation (prepending a prefix $v$ to any input string).

Suppose that $u = u_1 \ldots u_m$, $v = v_1 \ldots v_n$, where 
$u_1, \,\ldots\,, u_m, v_1, \,\ldots\,, v_n \in \{0,1\}$.
To express all functions of the form $[v \leftarrow u](.)$ over a finite
generating set, we observe that 

\smallskip

$[v \leftarrow u](.)$  $=$
$[v \leftarrow \varepsilon] \cdot [\varepsilon \leftarrow u](.)$,

\smallskip

$[\varepsilon \leftarrow u](.)$  $\,=\,$ 
$[\varepsilon \leftarrow u_m]$  $\cdot $  $[\varepsilon \leftarrow u_{m-1}]$
$\cdot $ \ $\dots$ \ $\cdot $ $[\varepsilon \leftarrow u_1](.)$,

\smallskip

$[v \leftarrow \varepsilon](.)$  $\,=\,$
$[v_1 \leftarrow \varepsilon]$  $\cdot$  $[v_2 \leftarrow \varepsilon]$  
$\cdot$  \ $\dots$ \ $\cdot$  $[v_n \leftarrow \varepsilon](.)$.

\smallskip

\noindent Thus all functions of the form $[v \leftarrow u](.)$ are expressed
over the set of four generators
$\,\{[0 \leftarrow \varepsilon],\, [1 \leftarrow \varepsilon],\,$
$[\varepsilon \leftarrow 0],\, [\varepsilon \leftarrow 1]\}$.

\smallskip

For all $w \in (\Gamma_M \cup \tau)^*$ and $x, y \in \{0,1\}^*$ 
we have:

\medskip

\hspace{0.5in}  $E_w(x) = y$ \ \ \ iff 
 \ \ \ $E_w \cdot {\sf id}_x(.) = [y \leftarrow x](.)$.

\medskip

\noindent So we have proved:

\begin{pro} \label{PROevMon} 
 \ The {\em evaluation problem} of the monoid $M_{2,1}$ over $\Gamma_M$ 
(or over $\Gamma_M \cup \tau$) reduces to the {\em word problem} of 
$M_{2,1}$ over $\Gamma_M$ (respectively over $\Gamma_M \cup \tau$) by a 
many-one $\log$-space reduction.
 \ \ \  \ \ \  $\Box$
\end{pro}

\subsection{Reduction of the evaluation problem to the word problem of 
$V$ over more general generating sets}

We saw already in Theorem \ref{redEVtoWP}(1) that the evaluation problem of 
$V$ over $\Gamma_{\!1} \cup \tau$ (and of $2V$ over $\Gamma_{\!2}$) reduces 
to the word problem of $V$ over $\Gamma_{\!1} \cup \tau$ (respectively $2V$ 
over $\Gamma_{\!2}$).
More generally, for $V$ we want to find a reduction from the evaluation 
problem to the word problem over a finite generating set $\Gamma_{\!1}$, or
over any other generating set containing $\Gamma_{\!1}$ (Theorem \ref{redEVtoWP}(2)).

Below, if $H$ is a subgroup of $V$ and $g \in V$, then $g \cdot H$ denotes
the left coset $\{g \cdot h :  h \in H\}$; and $H \cdot g$ is the right 
coset.

\begin{thm} \label{redMembersh} {\bf (commutation test for the evaluation
problem).} 

\smallskip

\noindent {\small \rm (1)} \ For any $g \in V$ and $x, y \in \{0,1\}^*$
we have:

\medskip

\hspace{0.3in} $g(x) = y$ \ \ \ iff 
 \ \ \ $g \cdot {\sf pFix}_V(x) \,=\, {\sf pFix}_V(y) \cdot g$.

\medskip

\noindent {\small \rm (2)} \ For any $x, y \in \{0,1\}^*$, let $\,\Gamma_x$,
$\Gamma_y$, $\Gamma_{\ov{x}}$, $\Gamma_{\ov{y}}$ be respectively 
generating sets of $\,{\sf pFix}_V(x)$, $\,{\sf pFix}_V(y)$, 
$\,{\sf pFix}_V(\ov{x})$, $\,{\sf pFix}_V(\ov{y})$. 
For any $g \in V$ and $x, y \in \{0,1\}^*$ we have:

\medskip

\hspace{0.3in} $g(x) = y$ \ \ \ iff     
 \ \ \ $(\forall \alpha \in \Gamma_x)(\forall \delta \in \Gamma_{\ov{y}})$
$[\,\delta \,g \,\alpha \,g^{-1}$  $=$  $g \, \alpha \, g^{-1} \,\delta\,]$ 
 \ \ {\rm and}

\smallskip

\hspace{1.45in}
$(\forall \beta \in \Gamma_y) (\forall \gamma \in \Gamma_{\ov{x}})$
$[\,\gamma \,g^{-1} \,\beta \,g$  $=$  $g^{-1} \,\beta \,g \,\gamma\,]$.
\end{thm}
{\sc Proof.} (1) Theorem \ref{redMembersh}(1) is equivalent to  

\smallskip

 \ \ \  $g(x) = y$ \ \ \ $\Leftrightarrow$
 \ \ \ ${\sf pFix}_V(x) \subseteq g^{-1} \cdot {\sf pFix}_V(y) \cdot g$
 \ and  
 \ ${\sf pFix}_V(y) \subseteq g \cdot {\sf pFix}_V(x) \cdot g^{-1}$.

\smallskip

\noindent The implication $[\Rightarrow]$ is straightforward, 

\smallskip

\noindent Let us prove $[\Leftarrow]$. 

For all $\alpha \in$ ${\sf pFix}(x)$: 
$ \ g\,\alpha\,g^{-1} \in {\sf pFix}(y)$.
Hence $\,g\,\alpha\,g^{-1}(y) = y$, hence 
$\alpha(g^{-1}(y)) = g^{-1}(y)$ for all $\alpha \in$ ${\sf pFix}(x)$;
therefore, ${\sf pFix}(x)$ $\subseteq$  ${\sf pFix}(g^{-1}(y))$.
By Lemma \ref{Fixes}, we conclude that $x \le_{\rm pref} g^{-1}(y)$. So,
$g^{-1}(y) = x u$ for some $u \in \{0,1\}^*$.

Similarly, for all $\beta \in {\sf pFix}(y)$: 
$ \ g^{-1}\,\beta\,g \,\in\, {\sf pFix}(x)$.
Therefore $\,g^{-1}\,\beta\,g(x) = x$, hence $\beta(g(x)) = g(x)$, so
${\sf pFix}(y)$  $\subseteq$  ${\sf pFix}(g(x))$.
Again by Lemma \ref{Fixes}, we conclude that $y \le_{\rm pref} g(x)$.
So, $g(x) = y v$ for some $v \in \{0,1\}^*$.

Since $g(x) = y v$ implies that $g(x)$ is defined, $g^{-1}(y) = x u$ now 
implies $y = g(x) \, u$. 
Now $g(x) = y v$ and $y = g(x) \, u$ imply that $g(x) = g(x) \, v \,u$.
Hence, $uv = \varepsilon$, hence $u = \varepsilon = v$.
So, $g(x) = y$.

\medskip

\noindent (2) The equality
$\,g \cdot {\sf pFix}_V(x) = {\sf pFix}_V(y) \cdot g\,$ is equivalent to the 
conjunction 

 \ \ \ $g \cdot {\sf pFix}_V(x) \cdot g^{-1} \subseteq {\sf pFix}_V(y)$ 
 \ and \ $g^{-1} \cdot {\sf pFix}_V(y) \cdot g \subseteq {\sf pFix}_V(x)$.

\noindent And for any $\varphi = \alpha_1 \,\ldots\, \alpha_m \in$ 
${\sf pFix}_V(x)$, with $\alpha_1, \,\ldots\,, \alpha_m \in \Gamma_x$ we 
have: 

 \ \ \ $g \varphi g^{-1}$  $=$ 
$g \alpha_1 g^{-1} \cdot g \alpha_2 g^{-1} \cdot$  $\,\ldots\,$ 
$\cdot g \alpha_m g^{-1}$.  

\noindent So, \ $(\forall \varphi \in {\sf pFix}_V(x))$ 
$[\,g\,\varphi\,g^{-1} \in {\sf pFix}_V(y)\,]$ 
 \ is implied by 
 \ $(\forall \alpha \in \Gamma_x)$ 
$[\,g\,\alpha\,g^{-1} \in {\sf pFix}_V(y)\,]$. 

\smallskip

\noindent Since $\Gamma_x \subseteq {\sf pFix}_V(x)$, the latter is also 
implied by the former. Hence:

$(\forall \varphi \in {\sf pFix}_V(x))$
$[\,g\,\varphi\,g^{-1} \in {\sf pFix}_V(y)\,]$
 \ is equivalent to 
 \ $(\forall \alpha \in \Gamma_x)$
$[\,g\,\alpha\,g^{-1} \in {\sf pFix}_V(y)\,]$.

\smallskip

\noindent The same reasoning applies to $\Gamma_y$ and
$ \ g^{-1} \cdot {\sf pFix}_V(y) \cdot g \subseteq {\sf pFix}_V(x)$. 
Hence:

$(\forall \varphi \in {\sf pFix}_V(y))$ 
$[\,g^{-1}\,\varphi\,g \in {\sf pFix}_V(x)\,]$  
 \ is equivalent to
 \ $(\forall \beta \in \Gamma_y):$
$[\,g^{-1}\,\beta\,g \in {\sf pFix}_V(x)\,]$.

\medskip 

To check membership of $g\,\alpha\,g^{-1}$ in ${\sf pFix}_V(y)$ 
we use the commutation test of Theorem \ref{CommTestP}(2). 
So, $g\,\alpha\,g^{-1} \in {\sf pFix}_V(y)$ is equivalent to
 \ $(\forall \delta \in \Gamma_{\ov{y}})$
$[\,\delta\,g\,\alpha \,g^{-1} = g \,\alpha\, g^{-1} \, \delta\,]$.

The same reasoning applies to 
$\,g^{-1}\,\beta\,g \subseteq {\sf pFix}_V(x)$.
 \ \ \  \ \ \  $\Box$

\bigskip

\noindent {\bf Proof of Theorem \ref{redEVtoWP}(2):}

\smallskip

\noindent Since $V$ is 2-generated (see \cite{DMason, BleakQuick}), we have
$ \ |\Gamma_x|$  $=$  $|\Gamma_y|$  $=$  $|\Gamma_{\ov{x}}|$  $=$
$|\Gamma_{\ov{y}}|$  $=$  $2$. 
And $\ov{x}, \ov{y}$ can be found from $x, y$ in log-space, by 
Prop.\ \ref{PROPexistcomplemSingl}.
Theorem \ref{redMembersh}(2) now gives eight equalities in $V$ whose 
conjunction holds iff $\,E_w(x) = y$.
Moreover, by Theorem \ref{FixP}(2) the elements of $\,\Gamma_x$,
$\Gamma_y$, $\Gamma_{\ov{x}}$, $\Gamma_{\ov{y}}$ can be computed from 
$x, y, \ov{x}, \ov{y}$ in log-space; more precisely, in log space one can
compute words over $\Gamma_{\!\!V}$ that represent the elements of 
$\,\Gamma_x$, $\Gamma_y$, $\Gamma_{\ov{x}}$, and $\Gamma_{\ov{y}}$.
Thus, there exists a log-space reduction from the evaluation problem
$E_w(x) = y$ (on input $w, x, y$) to the conjunction of eight word problems 
of $V$: 

\smallskip

$\delta_i\,w\,\alpha_j\,w^{-1} = w\,\alpha_j\,w^{-1}\,\delta_i$, 
 \ \ \ for $i,j \in\{1,2\}$,

\smallskip

$\gamma_i\,w^{-1}\,\beta_j\,w$  $=$  $w^{-1}\,\beta_j\,w\,\gamma_i$,
 \ \ \ for $i,j \in\{1,2\}$,

\smallskip

\noindent where $\alpha_1, \alpha_2 \in \Gamma_{\!\!V}^*$ are equal (in $V$) 
to the elements of $\Gamma_x$; $\,\beta_1, \beta_2 \in \Gamma_{\!\!V}^*$ are 
equal (in $V$) to the elements of $\Gamma_y$; 
$\,\gamma_1, \gamma_2 \in \Gamma_{\!\!V}^*$ are equal (in $V$) to the 
elements of $\Gamma_{\ov{x}}$; and $\delta_1, \delta_2 \in \Gamma_{\!\!V}^*$ 
are equal (in  $V$) to the elements of $\Gamma_{\ov{y}}$.

If $w \in \Gamma_{\!\!V}^*$ (or more generally, 
$w \in (\Gamma_{\!\!V} \cup \Delta)^*$), then this is a log-space eight-fold
conjunctive reduction of the evaluation problem over 
$\Gamma_{\!\!V} \cup \Delta$ to the word problem of $V$ over 
$\Gamma_{\!\!V} \cup \Delta$.  If $w \in (\Gamma_{\!\!V} \cup \Delta)^*$ then 
this is a log-space eight-fold conjunctive reduction to the word problem of 
$V$ over $\Gamma_{\!\!V} \cup \Delta$.
 \ \ \  \ \ \  $\Box$

\bigskip

\noindent The reduction of the evaluation problem to the 
word problem uses {\bf two commutation tests}:

\smallskip

\noindent $\bullet$ \ The first commutation test, for testing membership in 
a partial fixator (Theorem \ref{CommTestP}): 
 
 \ \ \  \ \ \ $g \,\in\, {\sf pFix}_V(P)$ \ \ iff \ \ $g$ commutes 
  with every element (or every generator) of $\,{\sf pFix}_V(Q)$.

\smallskip

\noindent $\bullet$ \ The second commutation test, for testing an evaluation 
relation (Theorem \ref{redMembersh}):

 \ \ \  \ \ \ $g(x) = y$ \ \ iff 
 \ \ $g \cdot {\sf pFix}_V(x) \,=\, {\sf pFix}_V(y) \cdot g$.

The latter can be reformulated in terms of the generators of the 
two partial fixators.

\section{Appendix}

\subsection{Evaluation problem for elements of $V$ and $2V$ given by 
tables.}

\smallskip

We can consider another version of the evaluation problem for $V$,
where every element $\varphi$ of $V$ are given by a table, instead of a 
word over a generating set of $V$.
A table is of the form $\{(u^{(i)}, v^{(i)}) : i = 1, \ldots, m\}$, which 
is a bijection between two finite maximal prefix codes 
$\{u^{(i)} : i = 1, \ldots, m\}$, $\{v^{(i)} : i = 1, \ldots, m\}$ 
$\subseteq \{0,1\}^*$.
The concepts of long and short data inputs also arises for tables, as
follows. 

A string $x \in \{0,1\}^*$ is a {\em long} data input for the above table
iff $u^{(i)} \le_{\rm pref} x$ for some $i = 1, \ldots, m$.
And $x$ is a {\em short} data input for the above table iff $x \in$
${\rm Dom}(\varphi)$, but $x$ is not long.

\medskip

\noindent The {\em evaluation problem of $V$ given by tables} is as follows.

\noindent {\sc Input:} \ A finite set of pairs 
$\{(u^{(i)}, v^{(i)}) : i = 1, \ldots, m\}$  $\subseteq$ 
$\{0,1\}^* \x \{0,1\}^*$, and $x, y \in \{0,1\}^*$.

\noindent {\sc Question 1 (general problem):} \ Is this set of pairs a table
for an element of $V$, and is $\,\varphi(x) = y$?

\noindent {\sc Question 2 (for long data inputs):}
 \ Is this set of pairs a table for an element of $V$, and is $x$ a long data 
input for that table, and is $\,\varphi(x) = y$ ?
Equivalently, the set of pairs a table, and do there exist $(u_j, v_j)$ in 
the table and $x_j, y_j \in \{0,1\}^*$ such that $x = u_j x_j$ and 
$y = v_j y_j$?

\medskip

\noindent The general evaluation problem with the empty string as data input
and output (i.e., $x = \e = y$) is equivalent to the {\em identity problem},
where the input is a table, and the question is whether that table represents 
the identity element of $V$. This problem is easy to solve, since
$\{(u^{(i)}, v^{(i)}) : i = 1, \ldots, m\}$ represents the identity iff
$u^{(i)} = v^{(i)}$ for $i = 1, \ldots, m$. This also requires checking 
whether $\{u^{(i)} : i = 1, \ldots, m\}$ is a maximal prefix code.
The general evaluation problem and the evaluation problem with long data 
inputs are certainly in {\sf P}.

\medskip

Similarly, an element of $2V$ could be given by a table that describes a 
bijection between finite maximal joinless codes. Just as for $V$, this 
leads to a general evaluation problem and the evaluation problem with 
long data inputs. These problems are easily seen to be in {\sf P}.

\subsection{The evaluation problem of the Thompson monoid $M_{2,1}$}

The evaluation problem and the concept of long and short data inputs, apply
to $M_{2,1}$ in the same way as for $G_{2,1}$.
We saw in Prop.\ \ref{PROevMon} that the evaluation problem for $M_{2,1}$ 
reduces to the word problem of $M_{2,1}$; this holds for finite generating 
sets $\Gamma_{\!M}$ and for circuit-like generating sets 
$\Gamma_{\!M} \cup \tau$. 

Similarly to Prop.\ \ref{VoverGammatoWP}(1), the word problem reduces to 
the evaluation problem with data input and output $\e$, for finite generating 
sets $\Gamma_{\!M}$ and for circuit-like generating sets 
$\Gamma_{\!M} \cup \tau$. We also have the analogue of Prop.\ 
\ref{VoverGammatoWP}(3) for $M_{2,1}$.

For long data inputs, the evaluation problem of $M_{2,1}$ is in {\sf DCF} 
over a finite generating set $\Gamma_{\!M}$, and it is {\sf P}-complete over
$\Gamma_{\!M} \cup \tau$.

For $M_{2,1}$ over $\Gamma_{\!M} \cup \tau$ the word problem is 
{\sf coNP}-complete. For $M_{2,1}$ over a finite generating set 
$\Gamma_{\!M}$, the word problem is in {\sf P} (see \cite{BiThompsMonV3}),
but it is still open whether it is {\sf P}-complete. 

All these problems can also be considered in the $2 M_{k,1}$, the monoid 
generalization of the Brin-Higman-Thompson group (see \cite{BinMk1}).

\subsection{Essential initial factor codes for defining ideal morphisms}

We usually define elements of $n V$ by tables that use finite maximal 
{\em joinless} codes, because such tables (when bijective) always define 
elements of $nV$.
On the other hand, for $n \ge 2$, tables based on finite essential initial 
factor codes do not always define elements of $nV$, as was proved in Prop.\ 
\ref{LEMuniqext}.
Nevertheless, since the unique maximum extension of a right ideal morphism 
typically requires a table with essential initial factor codes, we also need
to consider such tables.

Here we prove that one can decide in log-space whether a table based on
initial factor codes defines an element of $nV$. The concept of complementary
initial factor code is very useful for this, and will be introduced first.

\begin{defn} \label{DEFcomplEssInit}
 \ Two finite initial factor codes $P,Q \subseteq n\,A^*$ are {\em
complementary} \ iff 
 
\noindent {\small \rm (1)} \ $P\,(nA^*) \,\cap\, Q\,(nA^*) \,=\,\varnothing$,
 \ and 

\noindent {\small \rm (2)} \ $P \,\cup\, Q\,$ is an essential initial factor
code.
\end{defn}
If $P \subseteq n\,A^*$ is essential then $\varnothing$ is the unique 
complementary initial factor code of $P$. 
On the other hand, the complementary initial factor codes of $\varnothing$ 
are all the essential initial factor codes, e.g., $\{\e\}^n$. 

\bigskip

\noindent {\bf Notation:} 

\smallskip

$A_{\e} = \bigcup_{i=1}^n \{\e\}^{i-1} \x A \x \{\e\}^{n-i}$ \ \ \ (this is
the minimum generating set of the monoid $nA^*$).

\smallskip

\noindent For every set $Q \subseteq nA^*$, \ 

\smallskip

${\sf init}(Q)$  $\,=\,$ 
$\{z \in nA^* : (\exists q \in Q)[\, z \le_{\rm init} q\,]\,\}$
 \ \ \ (i.e., the set of initial factors of elements of $Q$); 

\smallskip

${\sf Sinit}(Q)$  $\,=\,$ ${\sf init}(Q) \minus Q$  \ \ \ (i.e., the set of
strict initial factors of elements of $Q$).

\begin{pro} \label{PROcomplEssInitp} {\bf (complementary initial factor
code).} 

\smallskip

\noindent Every finite initial factor code $P \subseteq n\,A^*$ has a 
(usually non-unique) complementary finite initial factor code 
$P' \subseteq n\,A^*$.  Moreover, $P'$ can be chosen so that

\smallskip

\noindent $\bullet$ \ \ $|P'| \ \le \ (|A|-1) \ \sum_{p \in P} |p|$, 

\smallskip

\noindent $\bullet$ \ \ ${\rm maxlen}(P') = {\rm maxlen}(P)$, \ and

\smallskip

\noindent $\bullet$ \ \ $P'$ can be computed from $P$ in log-space (if $P$
is given as a finite list of $n$-tuples of strings).
\end{pro}
{\sc Proof.} We define the set 

\medskip

 \ \ \   \ \ \ $P^{^{\#}} \ = \ $
$\big\{s \alpha \,:\, s \in {\sf Sinit}(P), \ \alpha \in A_{\e}, $
$ \ \{s \alpha\} \vee P = \varnothing \big\}$.

\medskip

\noindent The set $P^{^{\#}}$ is not necessarily an initial factor code (see 
the example after this proof), so we let 

\medskip

 \ \ \   \ \ \ $P' = \max_{\le_{\rm init}}(P^{^{\#}})$,

\medskip

\noindent i.e., $P'$ consists of the $\le_{\rm init}$-maximal elements of
$P^{^{\#}}$. 

We can compute $P'$ from $P^{^{\#}}$ in log-space by eliminating all 
elements of $P^{^{\#}}$ that have another element of $P^{^{\#}}$ as an 
initial factor.  In other words, $P'$ is the unique maximal initial factor 
code that is contained in $P^{^{\#}}$.
The formulas for $P^{^{\#}}$ and $P'$ immediately imply log-space 
computability of $P^{^{\#}}$ and $P'$, and the formulas for the cardinality 
$|P'|$ and for ${\rm maxlen}(P')$.

Let us check that $P'$ is a complementary initial factor code of $P$.
 
\smallskip

The condition $\,\{s \alpha\} \vee P = \varnothing\,$ in the formula for 
$P^{^{\#}}$ immediately implies $\,P\,nA^* \cap P^{^{\#}}\,nA^*$  $=$ 
$\varnothing$.  Since $P' \subseteq P^{^{\#}}$, we also have 
$\,P\,nA^* \cap P'\,nA^* = \varnothing$.

\medskip

We prove next that $P \cup P^{^{\#}}$ is essential. This follows if we show 
that for every $x = (x^{(1)}, \,\ldots\,, x^{(n)}) \in nA^{\ell}$ where 
$\ell = {\rm maxlen}(P) + 1$: 
 \ If $\,x \not\in P\,(nA^*)\,$ then $\,x \in P^{^{\#}}(nA^*)$.

Since $x \not\in P\,(nA^*)$, and $|x^{(h)}| = \ell > {\rm maxlen}(P)$ for 
all $h \in$ $\{1,\ldots,n\}$, it follows that for every $p \in P$ there 
exists $i \in \{1,\ldots,n\}$ such that we have: the longest common prefix 
of $x^{(i)}$ and $p^{(i)}$ is a {\sl strict} prefix of $p^{(i)}\,$ (in 
$A^*$).  This prefix is of course also a strict prefix of $x^{(i)}$, since 
$|x^{(i)}| = \ell$. We now consider a maximally long such prefix, i.e., a 
maximally long string in 

\smallskip

 \ \ \  \ \ \ $\{ r \in A^* :\,$ 
$i \in \{1,\ldots,n\},  \ p \in P, \ r \le_{\rm pref} x^{(i)}, \ $ 
$r < _{\rm pref} p^{(i)} \}$.

\smallskip

\noindent Let $j \in \{1,\ldots,n\}$ and $q \in P$ be a choice of values in 
$\{1,\ldots,n\}$ respectively $P$ where this maximum is reached. Then a 
maximally long $r$ is of the form $r = $
$\,x_1^{(j)} \,\ldots\, x_m^{(j)} < _{\rm pref} p^{(j)}$, 
$\,m = |r| < |p^{(j)}| < \ell$, where $\,x^{(j)} =$
$x_1^{(j)} \,\ldots\, x_m^{(j)} x_{m+1}^{(j)} \,\ldots\, x_{\ell}^{(j)}$.
Then $\,x_1^{(j)} \,\ldots\, x_m^{(j)} x_{m+1}^{(j)}\,$ is a prefix
of $x^{(j)}$ that (by maximality of the length $m$) is {\sl not} a prefix 
of $p^{(j)}$ for any $p \in P$.

For every $i \in \{1,\ldots,n\}$, let $s^{(i)}$ be the longest
common prefix of $x^{(i)}$ and $q^{(i)}$; for $i = j$ we already have
$s^{(j)} = x_1^{(j)} \,\ldots\, x_m^{(j)}$.
Then  \,$s = \big(s^{(1)}, \,\ldots\,, s^{(j-1)},$  
$\,x_1^{(j)} \,\ldots\, x_m^{(j)},$
$\,s^{(j+1)}, \,\ldots\,, s^{(n)}\big)$  $<_{\rm init} q$, hence  
$s \in {\sf Sinit}(P)$. Moreover, $\{s \alpha\} \vee P = \varnothing$,
where $\{\alpha\} = \{\e\}^{j-1} \x \{x_{m+1}^{(j)}\} \x \{\e\}^{n-j}$; 
this holds since $\,x_1^{(j)} \,\ldots\, x_m^{(j)} x_{m+1}^{(j)}\,$ 
is {\sl not} a prefix of $p^{(j)}$ for any $p \in P$.
Hence, $s \alpha \in P^{^{\#}}$, and $s \alpha \le_{\rm init} x$; thus, 
$x \in P^{^{\#}}(nA^*)$.  

From the fact that $P \vee P^{^{\#}}$ is essential it follows that 
$P \vee P'$ is essential, since every element of $P^{^{\#}}$ has a prefix in 
$P'$. I.e., $P^{^{\#}} \subseteq P'\,(nA^*)$, so 
$P^{^{\#}} \,(nA^*) \subseteq P'\,(nA^*)$. 
 \ \ \ $\Box$ 

\bigskip

\noindent {\bf Example where $P^{^{\#}}$ is not an initial factor code:}

\noindent Let $p = (11,00) \in 2A^*$ with $A = \{0,1\}$ and let $P = \{p\}$.
Let $s = (1,0)$ and $t = (1,00)$; so $s, t \in {\sf Sinit}(P)$. 
Let $\alpha = (0,\e) \in A_{\e}$.
Now $s \alpha = (10,0)$ and $t \alpha = (10,00)$.
Hence $\,s \alpha \vee p\,$ and $\,t \alpha \vee p\,$ do not exist, so
$\,s \alpha, t \alpha \in P^{^{\#}}$. But 
$\,s \alpha  <_{\rm init} t \alpha$; so $P^{^{\#}}$ is not an initial factor 
code in this example.

\begin{pro} \label{PROessDecid}
 \ It can be decided in log-space whether a finite initial factor code 
$P \subseteq nA^*$ (given as a finite list of $n$-tuples of strings) is 
{\em essential}.
\end{pro}
{\sc Proof.} By Prop.\ \ref{PROcomplEssInitp}, a complementary initial factor
code $P'$ of $P$ can be computed in log-space. And $P$ is essential iff 
$P' = \varnothing$. 
 \ \ \ $\Box$

\begin{pro} \label{PROfinitcodedec}
 \ Each one of the following questions can be decided in log-space.

\smallskip

\noindent {\sc Input:} A function $F: P \to Q$ from a finite initial factor 
code $P$ onto a finite initial factor code $Q$ (where $P$ and $Q$ are given 
as finite lists of $n$-tuples of strings).  Let $f: P\,(nA^*) \to Q\,(nA^*)$ 
be the right ideal ``morphism'' determined by the table $F$.  

\smallskip

\noindent {\sc Question 1:} Is $f$ a {\em function}?

\smallskip

\noindent {\sc Question 2:} Is $f$ {\em injective}?

\smallskip

\noindent {\sc Question 3:} Is $f$ {\em total}?

\smallskip

\noindent {\sc Question 4:} Is $f$ {\em surjective}?

\smallskip

\noindent {\sc Question 5:} Does $F$ define an element of the 
Brin-Higman-Thompson group $n G_{|A|,1}$?
\end{pro}
{\sc Proof.} (Q1) For every pair $p, p' \in P$ such that $p \ne p'$ and 
such that $\,p \vee p'\,$ exists, let $\,p \vee p' = pu = p'v$. 
If $f(p) \, u \ne f(p') \, v$ then $f$ is not a function.
If after checking all pairs $p.p'$ as above, no inequality was found, 
then $f$ is a function. \\        
(Q2) For this we check, as in (Q1), whether the inverse table of 
$F^{-1}: Q \to P$ is a function.  \\         
(Q3) We check whether $P$ is essential, using Prop.\ \ref{PROessDecid}. \\  
(Q4) We check whether $Q$ is essential, using Prop.\ \ref{PROessDecid}. \\  
(Q5) The table $F$ defines an element of $n G_{|A|,1}$ iff $F$ yields a 
{\sc yes} answer for all of the above questions.
 \ \ \ $\Box$

\bigskip

\bigskip




\bigskip

{\small
J.C.\ Birget

birget@camden.rutgers.edu
}


{\small

\begin{thebibliography}{99}

\bibitem{BiThomps} J.C.\ Birget, ``The groups of Richard Thompson and 
  complexity'', {\it International J.~of Algebra and Computation} 
  14(5, 6) (Dec.\ 2004) 569-626.  \ \ ( arxiv:math/0204292v2 )

\bibitem{BiCoNP} J.C.\ Birget, ``Circuits, coNP-completeness, and the groups
  of Richard Thompson'', {\it International J.~of Algebra and Computation}
  16(1) (Feb.\ 2006) 35-90. \ \ ( arxiv:math/0310335 )

\bibitem{BiThompsMon} J.C.~Birget, ``Monoid generalizations of the Richard
  Thompson groups'', {\it J.~of Pure and Applied Algebra} 213(2) (Feb.\
  2009) 264-278.

\bibitem{BiThompsMonV3} J.C.~Birget, ``Monoid generalizations of the Richard
  Thompson groups'', arxiv:0704.0189v3 \ (corrected version of 
  \cite{BiThompsMon}).

\bibitem{BinG} J.C.\ Birget, ``The word problem of the Brin-Thompson
  groups is coNP-complete'', {\it J.\ Algebra} 553 (1 July 2020) 268-318. 
  \ \ ( arxiv.org/1902.03852 )

\bibitem{BinGk1} J.C.\ Birget, ``The word problem of the Brin-Higman-Thompson
  groups'', arxiv.org/2006.14968  

\bibitem{BinMk1} J.C.\ Birget, ``A monoid version of the Brin-Higman-Thompson
  groups'', arxiv.org/2006.15355 

\bibitem{BleakQuick} C.\ Bleak, M.\ Quick, ``The infinite simple group $V$
  of Richard J.\ Thompson: presentations by permutations'',
  {\it Groups, Geometry, and Dynamics} 11 (2017) 1401-1436.

\bibitem{CFP} J.W.\ Cannon, W.J.\ Floyd, W.R.\ Parry,
  ``Introductory notes on Richard Thompson's groups'', \\
  {\it L'Enseignement Math\'ematique} 42 (1996) 215-256.

\bibitem{Cook} S.A.\ Cook, ``The complexity of theorem-proving procedures'',
  Proceedings 3rd ACM STOC (1971) 151–158.

\bibitem{DuKo} D.Z.\ Du, K.I.\ Ko, {\it Theory of Computational 
  Complexity}, Wiley (2000).

\bibitem{GHR} R.\ Greenlaw, H.J.\ Hoover, W.L.\ Ruzzo, {\it Limits to 
  parallel computation: {\rm P}-completeness theory}, Oxford U.P.\
  (1995). 

\bibitem{Harrison} M.A.\ Harrison, {\it Introduction to formal language 
  theory}, Addison-Wesley (1978). 

\bibitem{HemaOgi} L.H.\ Hemaspaandra, M.\ Ogihara, {\it The complexity
  theory companion}, Springer (2002).

\bibitem{Hig74} G.\ Higman, ``Finitely presented infinite simple groups'',
  Notes on Pure Mathematics 8, The Australian National University,
  Canberra (1974).

\bibitem{HU} J.E.\ Hopcroft, J.D.\ Ullman, {\it Introduction to
  automata theory, languages, and computation}, Addison-Wesley
  (1979 edition).

\bibitem{Ladner} R.E.\ Ladner, ``The circuit value problem is log space 
  complete for P'',  {\it ACM Sigact News} 7(1) (1975) 18-20.

\bibitem{LS} J.\ Lehnert, P.\ Schweitzer, ``The co-word problem for 
  the Higman-Thompson group is context-free'', 
  {\it Bull.\ London Math.\ Soc.} 39 (2007) 235–241.  

\bibitem{DMason} D.R.\ Mason, ``On the 2-generation of certain finitely
  presented infinite simple groups'', {\it Journal of the London
  Mathematical Society \ s2}-16(2), (Oct.\ 1977) 229-231.

\bibitem{Papadim} Ch.\ Papadimitriou, {\it Computational complexity},
  Addison-Wesley (1994).

\end{thebibliography}

}
\end{document}